\documentclass[12pt]{article}
\newcommand{\bR}{{\bf R}}
\newcommand{\bN}{{\bf N}}
\newcommand{\bZ}{{\bf Z}}
\newcommand{\bC}{{\bf C}}
\newcommand{\vro}{{\varrho}}
\newcommand{\eps}{{\varepsilon}}
\newcommand{\iy}{{\infty}}

\newcommand{\al}{{\alpha}}
\newcommand{\be}{{\beta}}

\newcommand{\la}{{\lambda}}

\newcommand{\cB}{{\cal B}}

\newcommand{\cM}{{\cal M}}

\newcommand{\cR}{{\cal R}}

\newtheorem{lemma}{Lemma}[section]
\newtheorem{proposition}[lemma]{Proposition}
\newtheorem{theorem}[lemma]{Theorem}
\newtheorem{corollary}[lemma]{Corollary}

\title
{\normalsize\bf
\vskip 2truecm
ONE-SIDED INVERTIBILITY OF BINOMIAL FUNCTIONAL OPERATORS\\ WITH A
SHIFT IN REARRANGEMENT-INVARIANT SPACES }
\author
{\normalsize
ALEXEI YU. KARLOVICH\footnote{Partially supported by F.C.T.
(Portugal) grant PRAXIS XXI/BPD/22006/99.} \ and YURI
I.~KARLOVICH\footnote{Partially supported by CONCACYT (M\'exico)
grant, C\'atedra Patrimonial, No. 990017-EX., nivel II.} }

\date{}
\begin{document}
\maketitle
\footnotesize
\baselineskip=12pt
\begin{quote}
Let $\Gamma$ be an oriented Jordan smooth curve and $\alpha$ a
diffeomorphism of $\Gamma$ onto itself which has an arbitrary
nonempty set of periodic points. We prove criteria for 
one-sided invertibility of the binomial functional operator
\[
A=aI-bW
\]
where $a$ and $b$ are continuous functions, $I$ is the identity
operator, $W$ is the shift operator, $Wf=f\circ\alpha$, in a
reflexive rearrangement-invariant space $X(\Gamma)$ with Boyd
indices $\alpha_X,\beta_X$ and Zippin indices $p_X,q_X$ satisfying
inequalities
\[
0<\alpha_X=p_X\le q_X=\beta_X<1.
\]
\end{quote}
\normalsize

\vskip 1truecm
\baselineskip=15pt
\section{Introduction}
Let $\Gamma$ be an oriented Jordan (i.e., homeomorphic to a
circle) smooth curve. Let $\alpha$ be a diffeomorphism of $\Gamma$
onto itself which preserves or changes the orientation on
$\Gamma$. We consider the binomial functional operator
\begin{equation}\label{eq:operator}
A:=aI-bW
\end{equation}
acting in a rearrangement-invariant space $X(\Gamma)$, where
$a$ and $b$ are continuous functions on $\Gamma$, $I$ is
the identity operator, and $W$ is the shift operator defined by
\[
(Wf)(t):=f[\alpha(t)],\quad t\in\Gamma.
\]

An investigation of the two- and one-sided invertibility of
functional operators (in particular, (\ref{eq:operator})) in
various functional spaces plays an important role in the theory of
functional differential operators (see, e.g., \cite{Antonevich},
\cite{AntLeb}, \cite{Kurbatov}), theory of singular integral
operators, convolution type operators and pseudodifferential
operators with shifts and/or oscillating coefficients (see
\cite{AntLebBel}, \cite{KrLit}, \cite{Lit} and the
references therein), theory of dynamical systems \cite{Lat},
theory of Banach lattices and Banach $C(K)$-modules \cite{AAK}, etc.

Criteria for the two-sided invertibility of the operator
(\ref{eq:operator}) in Lebesgue spaces $L^p(\Gamma), 1<p<\infty$, 
were obtained by V.~G.~Kravchenko and the second author
\cite{KK76}. Criteria for one-sided invertiblity of
(\ref{eq:operator}) in $L^p(\Gamma)$ were established by R.~Mardiev
\cite{Mar85,Mar88}. These results were extended by V.~Aslanov and
the second author to the case of reflexive Orlicz spaces
$L^M(\Gamma)$ and were announced in \cite{AK89},
the full proofs were given in \cite[Section~4]{AK2000}. 
All these results are contained in the survey \cite{Karl95} 
of the second author.
The paper \cite{AKarl00} of the
first author was devoted to a further generalization and
development of ideas and results of \cite{AK2000,AK89} to the case
of reflexive rearrangement-invariant spaces and shifts having only
two fixed points.

In this paper we extend the results of \cite{AKarl00} to the case
of shifts $\alpha$ preserving or changing the orientation on $\Gamma$
and having an arbitrary nonempty set of periodic points. A nontrivial
example of such shift was constructed in \cite{KK76} (see also
\cite[p. 74]{KrLit}). We obtain criteria for one-sided
invertibility of the functional operator (\ref{eq:operator}) in a
reflexive rearrangement-invariant space $X(\Gamma)$ of fundamental
type with nontrivial Boyd indices. Such spaces are wide
generalizations of Lebesgue and Orlicz spaces. 

The paper is organized as follows. In Section 2 we formulate
necessary properties of rearrangement-invariant spaces and their
interpolation characteristics (the Boyd and Zippin indices). We
also describe the structure of the set $\Lambda$ of periodic
points of the shift $\alpha$ (which have the same multiplicity
$m\in{\bf N}$ if $\alpha$ preserves the orientation on $\Gamma$ and
have multiplicities 1 and 2 otherwise, in the latter case we put
$m=2$). Difficulties appear in the case of the infinite boundary
$\partial\Lambda$ of the set of periodic points. In that case the
curve $\Gamma$ to within a finite subset can be represented as a
finite union of pairwise disjoint open arcs of three types: on the
arcs of first type the shift $\alpha$ is Carleman (that is,
$\alpha_m(t)\equiv t$), on the arcs of second type the shift
$\alpha_m$ has only two fixed points (the endpoints of the arc),
and the closure of the union of the arcs of third type contains
the set $(\partial\Lambda)\,'$ of all limit points of the boundary
$\partial\Lambda$ of $\Lambda$. Moreover, we can choose the arcs
of third type so small as we want. All these arcs are invariant
with respect to the shift $\alpha_m$.

In Section 3 we prove sufficient conditions for the two-sided
invertibility of the functional operator (\ref{eq:operator}) in
$X(\Gamma)$. These results are based on an estimate from above for
the spectral radius of the weighted shift operator $gW$ ($g\in
C(\Gamma)$) in rearrangement-invariant spaces, which is obtained
with the help of interpolation from known results for Lebesgue
spaces.

Section 4 is devoted to criteria for one-sided invertibility
of the functional operator (\ref{eq:operator}) in $X(\Gamma)$ in
the case of an arbitrary nonempty set of fixed points of $\alpha$.
In their proofs we essentially use a decomposition of $\Gamma$
into a union of arcs of three types. Further, for the functional
operator (\ref{eq:operator}) we define five sets $\Gamma_j , \:
j\in\{1,2,\ldots,5\}$, and in terms of these sets we prove a
criterion for one-sided invertibility of $A$. Roughly speaking,
$\Gamma_1$ is a Carleman part of $\Gamma$, the sets $\Gamma_2$ and
$\Gamma_3$ control the two-sided invertibility of $A$ if,
respectively, the coefficient $a$ or $b$ dominates, the set
$\Gamma_4$ controls the right invertibility of $A$, and the set
$\Gamma_5$ controls the left invertibility of $A$. Stress also
that according to Corollary \ref{le:onesided-spectrum}, if
$\alpha$ has a finite set of fixed points on $\Gamma$, then for
considered rearrangement-invariant spaces (with the non-coinciding Boyd indices)
in contrast to Lebesgue
spaces, the intersection of the left and right spectra of the
shift operator $W$ consists of a finite union of annuli depending
on the values of the derivative $\alpha'$ at the fixed points of
$\alpha$ and the both Boyd indices. In the case of
Lebesgue spaces these annuli degenerate into circles. This shows a
new quality arising for functional operators in
rearrangement-invariant spaces in view of the non-coincidence, in general, 
of the Boyd indices.

In Section 5 we extend the results of Section 4 to the case of
periodic points of arbitrary multiplicity $m$. The proof is based
on the equivalence of the right invertibility of $A$ and the
conjunction of the right invertibility of an operator of the form
(\ref{eq:operator}), with the shift $\alpha_m$ having only fixed
points, and some additional conditions. Using that equivalence, we
reduce the general case to the case of fixed points. The left
invertiblity is studied by passing to adjoint operators and making
use of the reflexivity of the space $X(\Gamma)$. At the end of
this section we calculate the spectrum of the weighted shift
operator $gW$ with the continuous coefficient $g$. 
\section{Spaces and shifts}
\setcounter{equation}{0}
\subsection{Rearrangement-invariant spaces}
For a general discussion of rearrangement-invariant spaces, see
\cite{BeSh,kps,LT}. In this subsection we collect
necessary facts in the abstract setting of (finite) measure spaces.

Let $(\cR,\mu)$ be a nonatomic finite measure space.
Denote by $\cM=\cM(\cR,\mu)$ the
set of all $\mu$-measurable complex-valued functions on $\cR$, and let
$\cM^+$ be the subset of functions from $\cM$ whose values lie in $[0,\iy]$.
The characteristic function of a $\mu$-measurable set $E\subset\cR$ will be
denoted by $\chi_E$.
A mapping $\rho:\cM^+\to [0,\iy]$ is called a function norm if for all
functions $f,g,f_n \in \cM^+ \ (n\in{\bf N})$, for all constants $a\ge 0$ and for
all $\mu$-measurable subsets $E$ of $\cR$, the following properties hold:
\begin{eqnarray*}
(i) & &\rho(f)=0 \ \Leftrightarrow \ f=0 \ \mu\mbox{-a.e.},
\       \rho(af)=a\rho(f),
\        \rho(f+g) \le \rho(f)+\rho(g),\\
(ii) & &0\le g \le f \ \mu\mbox{-a.e.} \ \Rightarrow \ \rho(g) \le \rho(f)
\quad\mbox{(the lattice property)},\\
(iii) & &0\le f_n \uparrow f \ \mu\mbox{-a.e.} \ \Rightarrow \
       \rho(f_n) \uparrow \rho(f)\quad\mbox{(the Fatou property)},\\
(iv) & &\rho(\chi_E) <\iy,\quad \int_E f\,d\mu \le C_E\rho(f)
\end{eqnarray*}
with $C_E \in (0,\iy)$ depending on $E$ and $\rho$ but independent of $f$.
The collection $X=X(\rho)$ of all functions $f\in\cM$
for which $\rho(|f|)<\iy$ is called a Banach function space. For each
$f \in X$, the norm of $f$ is defined by
\[
\|f\|_X :=\rho(|f|).
\]

If $\rho$ is a function norm, its associate norm $\rho'$ is
defined on $\cM^+$ by
\[
\rho'(g):=\sup\left\{
\int_\cR fg\,d\mu \ : \ f\in \cM^+, \ \rho(f) \le 1
\right\}, \quad g\in \cM^+.
\]
The Banach function space $X(\rho')$ determined by the function norm
$\rho'$ is called the associate space of $X=X(\rho)$ and is denoted by $X'$.
The associate space $X'$ is a subspace of the dual space $X^*$.

In the following we will consider only separable measure spaces.
Note that the Lebesgue measure is separable  (for the definition 
and the proof of this fact, see, e.g., \cite[Ch.~1, Subsection~6.10]{ka}).
\begin{lemma}$\!\!\!${\bf .}\label{le:ass&adj}
Let $\mu$ be a separable measure.

{\rm (a)} A Banach function space $X$ is separable if and only if its associate
space $X'$ is canonically isometrically isomorphic to the dual space $X^*$ of
$X$.

{\rm (b)} A Banach function space $X$ is reflexive if and only if both
$X$ and its associate space $X'$ are separable.
\end{lemma}

This lemma follows from Corollaries~4.3, 4.4 and 5.6 \cite[Ch.~1]{BeSh}.

Let $\cM_0$ and $\cM_0^+$ be the classes of $\mu$-a.e. finite
functions in $\cM$ and $\cM^+$, respectively. Two functions
$f,g\in\cM_0$ are said to be equimeasurable if
\[
\mu\{x\in\cR\::\:|f(x)|>\la\}=\mu\{x\in\cR\::\:|g(x)|>\la\}
\quad\mbox{for all}\quad\la\ge 0.
\]
A function norm $\rho:\cM^+ \to [0,\iy]$ is said to be
rearrangement-invariant, if $\rho(f)=\rho(g)$ for every pair of equimeasurable
functions $f,g \in \cM^+_0$. In that case, the Banach function space
$X=X(\rho)$ generated by $\rho$ is said to be the rearrangement-invariant
space (briefly r.-i. space). Lebesgue, Orlicz, Lorentz spaces are important
classical examples of r.-i. spaces.
\subsection{Boyd and Zippin indices}
A measurable function $\vro:(0,\iy)\to(0,\iy)$ is said to be
submultiplicative if
\[
\vro(x_1x_2)\le\vro(x_1)\vro(x_2)
\quad\mbox{for all}\quad x_1,x_2\in(0,\iy).
\]
The behavior of the measurable submultiplicative function $\vro$
in neighborhoods of zero and infinity is described by the
quantities (see \cite[Ch.~2, Theorem~1.3]{kps})
\begin{equation}\label{eq:indices-def}
\al(\vro):=\sup_{x\in(0,1)}\frac{\log\vro(x)}{\log x}
=\lim_{x\to 0}\frac{\log\vro(x)}{\log x},
\quad
\be(\vro):=\inf_{x\in(1,\iy)}\frac{\log\vro(x)}{\log x}
=\lim_{x\to \infty}\frac{\log\vro(x)}{\log x}.
\end{equation}
One can prove that $\al(\vro) \le \be(\vro)$ and these numbers are finite.
The numbers $\al(\vro)$ and $\be(\vro)$ are called the lower and upper
indices of the measurable submultiplicative function $\vro$.

The idea of using indices of some submultiplicative functions for the
description of properties of Orlicz spaces goes back to W.~Matuszewska
and W.~Orlicz, 1960. Matuszewska-Orlicz indices were generalized by
D.~W.~Boyd and M.~Zippin to the case of rearrangement-invariant spaces
(for the history and precise references, see \cite{mal}).

By the Luxemburg representation theorem \cite[Ch.~2, Theorem~4.10]{BeSh},
there is the unique rearrangement-invariant function norm
$\overline{\rho}$ over $[0,\mu(\cR)]$ with Lebesgue measure $m$ such that
\[
\rho(f) = \overline{\rho}(f^*) \quad\mbox{for all}\quad f\in \cM_0^+
\]
where $f^*$ is the non-increasing rearrangement of $f$
(see, e.g., \cite[p.~39]{BeSh}). The r.-i. space over $([0,\mu(\cR)],m)$
generated by $\overline{\rho}$ is called the Luxemburg representation
of $X$ and is denoted by $\overline{X}$. For each $x>0$, let $E_x$
denote the dilation operator defined on $\cM_0([0,\mu(\cR)],m)$ by
\begin{equation}\label{eq:dilation}
(E_x f)(t):=
\left\{
\begin{array}{ll}
f(xt), & xt\in[0,\mu(\cR)]\\
0,     & xt\not\in[0,\mu(\cR)]
\end{array}
\right.
, \quad t\in [0,\mu(\cR)].
\end{equation}
Together with $E_x$ consider the family of operators
\[
E_x^{(\lambda)}:=\Pi_\lambda E_x\Pi_\lambda,
\quad
\lambda\in(0,\mu(\cR)],
\]
where $\Pi_\lambda:=\chi_\lambda I$ and $\chi_\lambda$ is the
characteristic function of the segment $[0,\lambda]$. Consider
the function
\[
h_X(x,\lambda):=\|E_{1/x}^{(\lambda)}\|_{{\cal B}(\overline{X})},
\quad x\in(0,\infty),\quad \lambda\in(0,\mu(\cR)],
\]
where ${\cal B}(\overline{X})$ is the Banach algebra of the bounded
linear operators on $\overline{X}$.
\begin{lemma}$\!\!\!${\bf .}\label{le:properties}
{\rm (a)}
For every $x\in(0,\infty)$, we have
$h_X(x,\mu(\cR))\le\max\{1,x\}$;

\noindent
{\rm (b)}
for $x \in (0, \infty)$, the function $h_X(x,\lambda)$ is non-decreasing
in $\lambda\in(0,\mu(\cR)]$;

\noindent {\rm (c)} for $\lambda \in (0, \mu(\cR)]$, the function
$h_X(x,\lambda)$ is non-decreasing and submultiplicative in \\
$x\in(0,\infty)$;

\noindent
{\rm (d)}
if $X'$ denotes the associate space of $X$, then
\begin{equation}\label{eq:properties-1}
h_X(x,\lambda)=xh_{X'}\left(\frac{1}{x},\lambda\right),
\quad x\in(0,\infty),\quad \lambda\in(0,\mu(\cR)];
\end{equation}

\noindent {\rm (e)} if $x\in(0,1]$ and $\lambda\in(0,\mu(\cR)/2]$,
then $h_X(x,2\lambda)\le 2h_X(x,\lambda)$;

\noindent
{\rm (f)}
if $0<\lambda\le\nu\le\mu(\cR)$, then
there is a constant $C_{\lambda,\nu}>0$
such that
\[
h_X(x,\lambda)
\le
h_X(x,\nu)
\le
C_{\lambda,\nu} h_X(x,\lambda),
\quad x\in(0,\infty).
\]
\end{lemma}
{\bf Proof.}
The statement (a) is well-known (see, e.g., \cite[p.~165]{BeSh}).

(b) If $\lambda_1<\lambda_2$, then
$\Pi_{\lambda_1}=\Pi_{\lambda_1}\Pi_{\lambda_2}=\Pi_{\lambda_2}\Pi_{\lambda_1}$.
Hence,
\begin{equation}\label{eq:properties-2}
E_{1/x}^{(\lambda_1)}
=
\Pi_{\lambda_1}E_{1/x}\Pi_{\lambda_1}
=
\Pi_{\lambda_1}\Pi_{\lambda_2}
E_{1/x}
\Pi_{\lambda_2}\Pi_{\lambda_1}
=
\Pi_{\lambda_1}
E_{1/x}^{(\lambda_2)}
\Pi_{\lambda_1}.
\end{equation}
Since $\|\Pi_{\lambda_1}\|_{\cB(\overline{X})}\le 1$, we infer
from (\ref{eq:properties-2}) that
\[
h_X(x,\lambda_1)
=
\|E_{1/x}^{(\lambda_1)}\|_{\cB(\overline{X})}
\le
\|E_{1/x}^{(\lambda_2)}\|_{\cB(\overline{X})}
=
h_X(x,\lambda_2).
\]

(c) and (d). For every $x\in(0,\infty)$ and
every $\lambda\in(0,\mu(\cR)]$, we have
\begin{equation}\label{eq:properties-3}
(E_x^{(\lambda)}f)^*(t)\le E_x^{(\lambda)} f^*(t),
\quad t\in(0,\mu(\cR)].
\end{equation}
Using this inequality, one can obtain by analogy with the remark after
\cite[Ch.~3, Corollary~6.11]{BeSh} that the function $h_X(x,\lambda)$
is non-decreasing and submultiplicative in $x\in(0,\infty)$. Moreover,
(\ref{eq:properties-1}) holds.

(e) From (\ref{eq:properties-2}), (\ref{eq:properties-3}) and the
monotonicity of $f^*$ we see that for $t\in (0,\mu(\cR)]$,
\begin{eqnarray}
&&
(E_{1/x}^{(2\lambda)}f)^*(t)
\le
E_{1/x}^{(2\lambda)}f^*(t)
=
E_{1/x}^{(2\lambda)}(\chi_\lambda f^*)(t)+
E_{1/x}^{(2\lambda)}(\chi_{[\lambda,2\lambda]} f^*)(t)
\nonumber\\
&&=
\chi_\lambda(t)E_{1/x}^{(2\lambda)}(\chi_\lambda f^*)(t)+
E_{1/x}^{(2\lambda)}(\chi_{[\lambda,2\lambda]} f^*)(t)
=
E_{1/x}^{(\lambda)}f^*(t)+
E_{1/x}^{(2\lambda)}(\chi_{[\lambda,2\lambda]} f^*)(t),
\label{eq:properties-4}
\\[3ex]
&&
\Big(E_{1/x}^{(2\lambda)}(\chi_{[\lambda,2\lambda]} f^*)\Big)^*(t)
\le
E_{1/x}^{(2\lambda)}(\chi_{[\lambda,2\lambda]} f^*)^*(t)
=
E_{1/x}^{(2\lambda)}\Big(\chi_\lambda(t) f^*(t+\lambda)\Big)
\nonumber\\
&&
\le
E_{1/x}^{(2\lambda)}(\chi_\lambda f^*)(t)
=
\chi_\lambda(t)E_{1/x}^{(2\lambda)}(\chi_\lambda f^*)(t)
=
E_{1/x}^{(\lambda)}f^*(t).
\label{eq:properties-5}
\end{eqnarray}
Inequalities (\ref{eq:properties-4}), (\ref{eq:properties-5}) and
properties of the r.-i. invariant function norm $\overline{\rho}$ give
\begin{equation}\label{eq:properties-6}
\overline{\rho}(|E_{1/x}^{(2\lambda)}f|)
=
\overline{\rho}((E_{1/x}^{(2\lambda)}f)^*)
\le
2 \overline{\rho}(E_{1/x}^{(\lambda)}f^*).
\end{equation}
Obviously,
$
\{f^*:f \in \overline{X},\, \overline{\rho}(|f|)\le 1\}\subset
\{|f|:f \in \overline{X},\, \overline{\rho}(|f|)\le 1\}.
$
Hence, from (\ref{eq:properties-6}) we get
\begin{eqnarray*}
h_X(x,2\lambda)= \|E_{1/x}^{(2\lambda)}\|_{\cB(\overline{X})}
=
\sup_{f\in\overline{X},\, \overline{\rho}(|f|)\le 1}
\overline{\rho}(|E_{1/x}^{(2\lambda)}f|)
\le
2\sup_{f\in\overline{X},\, \overline{\rho}(|f|)\le 1}
\overline{\rho}(E_{1/x}^{(\lambda)}f^*)\qquad
\\
 \le2\sup_{f\in\overline{X},\, \overline{\rho}(|f|)\le 1}
\overline{\rho}(E_{1/x}^{(\lambda)}|f|)
=2\sup_{f\in\overline{X},\, \overline{\rho}(|f|)\le 1}
\overline{\rho}(|E_{1/x}^{(\lambda)}f|)
=
2\|E_{1/x}^{(\lambda)}\|_{\cB(\overline{X})}
=
2h_X(x,\lambda).
\end{eqnarray*}

(f) Choose $n\in\bN$ such that $\nu/2^n\le\lambda<\nu/2^{n-1}$.
Then from (b) and (e) we get
\begin{equation}\label{eq:properties-7}
h_X(x,\lambda) \le h_X(x,\nu) \le 2^nh_X(x,\nu/2^n) \le
2^nh_X(x,\lambda), \quad x\in(0,1].
\end{equation}
If $x\in(1,\infty)$, then from (\ref{eq:properties-1}) and
(\ref{eq:properties-7}) we get the same inequality
(\ref{eq:properties-7}) for $x\in(1,\infty)$.
\rule{2mm}{2mm}

The indices of the non-decreasing (and hence, measurable) and
submultiplicative function $h_X(\cdot,\mu(\cR))$ are called the
Boyd indices of the r.-i. space $X$ \cite{b4} and denoted by
\[
\alpha_X:=\alpha(h_X(\cdot,\mu(\cR))), \quad
\beta_X:=\beta(h_X(\cdot,\mu(\cR))).
\]
From Lemma~\ref{le:properties}(a),~(d) and equalities
(\ref{eq:indices-def}) it follows that
\[
0\le\alpha_X\le\beta_X\le 1,
\]
\begin{equation}\label{eq:ind-ass}
\alpha_X+\beta_{X'}=\alpha_{X'}+\beta_X=1.
\end{equation}

For each $t\in(0,\mu(\cR)]$, let $\Omega$
be a $\mu$-measurable subset of $\cR$ with $\mu(\Omega)=t$, and let
$\varphi_X(t):=\|\chi_\Omega\|_X$. The function $\varphi_X$ so defined
is called the fundamental function of the r.-i. space $X$. Put
\[
M_X(x):=\limsup_{t\to 0} \frac{\varphi_X(xt)}{\varphi_X(t)},
\quad x\in(0,\iy).
\]
The function $M_X$ is non-decreasing (and hence, measurable) and
submultiplicative (see \cite[Section~4]{mal}). The indices of this
function are called the Zippin (or fundamental) indices of the
r.-i. space $X$ \cite{zippin} and denoted by $p_X:= \al(M_X),q_X:=
\be(M_X)$. Generally, it can be proved (see, e.g.,
\cite[Section~4]{mal}) that 
\[
\al_X \le p_X \le q_X \le \be_X.
\]
\begin{lemma}$\!\!\!${\bf .}\label{le:subset}
Let $({\cal R},\mu)$ be a finite measure space and $\Omega$ be a
$\mu$-measurable subset of ${\cal R}$ with
$\lambda:=\mu(\Omega)>0$. Let $X$ be an r.-i. space over $({\cal
R},\mu)$ generated by an r.-i. function norm $\rho$. Then its
subspace $P_\Omega X$, where $P_\Omega:=\chi_\Omega I$, is an
r.-i. space over $(\Omega,\mu)$ generated by the same r.-i.
function norm $\rho$, but defined on functions with support in
$\Omega$. Moreover, $X$ and $P_\Omega X$ have the same Boyd and
Zippin indices.
\end{lemma}
{\bf Proof.} It is easy to check that an r.-i. function norm
$\rho$ defined on $\cM^+(\cR,\mu)$ is a Banach function norm on
$\cM^+(\Omega,\mu)$ as well. Moreover, $\rho$ is an r.-i. function
norm on $\cM^+(\Omega,\mu)$. Hence, $P_\Omega X$ is an r.-i.
space.

It is easy also to see that $\Pi_\lambda\overline{X}$ is the
Luxemburg representation for the subspace $P_\Omega X$, and
\begin{equation}\label{eq:subset-1}
\|E_{1/x}^{(\lambda)}\|_{\cB(\overline{X})}
=
\|E_{1/x}\|_{\cB(\Pi_\lambda\overline{X})},
\end{equation}
where the operator $E_{1/x}$ from the right is defined by
(\ref{eq:dilation}) on functions in $\cM_0^+([0,\lambda],m)$. From
(\ref{eq:indices-def}), (\ref{eq:subset-1}) and
Lemma~\ref{le:properties}(f) it follows that the Boyd indices of
$X$ and $P_\Omega X$ coincide.

On the other hand, $\varphi_X(t)=\varphi_{P_\Omega X}(t)$ for every
$t\in(0,\lambda]$. Hence, $M_X(x)=M_{P_\Omega X}(x)$ for every
$x\in (0,\infty)$. Consequently, the Zippin indices of $X$ and $P_\Omega X$ coincide too.
\rule{2mm}{2mm}

An r.-i. space $X$ is said to be of fundamental type if its Boyd and Zippin
indices coincide:
\[
\al_X=p_X,\quad\quad\be_X=q_X.
\]
Lebesgue, Orlicz, and Lorentz spaces are examples of spaces of fundamental
type \cite{feh83}. For the Lebesgue spaces $L^p, 1\le p\le \iy$,
all indices are equal to $1/p$. But there are r.-i. spaces for which the
Boyd and Zippin indices do not coincide, that is, there exist r.-i. spaces
of non-fundamental type (see \cite{mal} and the references given there).
We will say that the Boyd indices are nontrivial if
\[
0<\al_X\le\be_X<1.
\]
In the case of Orlicz spaces these inequalities are equivalent to
the reflexivity of the space (see, e.g., \cite{mal}). Examples of
Young functions which generate reflexive Orlicz spaces whose Boyd
indices do not coincide are given in \cite{AK89} and \cite[p.~93]{mal1}.

Boyd indices are interpolation characteristics of r.-i. spaces.
\begin{theorem}\label{th:spectral-est}
{\rm (see \cite{b4} and \cite[Theorem~2.2]{AKarl00})}. Let $X$ be
an r.-i. space with nontrivial Boyd indices. If there are numbers
$p_0,p_1\in(1,\infty)$ such that $$ 1/p_1<\alpha_X\le
\beta_X<1/p_0, $$ and if an operator $T$ is bounded in the
Lebesgue spaces $L^{p_0}$ and $L^{p_1}$, then $T$ is bounded in
$X$ and
\begin{equation}\label{eq:spectral-est}
r(T;X)\le\max\Big\{r(T;L^{1/\alpha_X}),r(T;L^{1/\beta_X})\Big\},
\end{equation}
where $r(T;X)$ is the spectral radius of $\, T$ in the r.-i. space
$X$, and $r(T;L^{1/\alpha_X})$ and $r(T;L^{1/\beta_X})$ are the
spectral radii of $\; T$ in the Lebesgue spaces $L^{1/\alpha_X}$
and $L^{1/\beta_X}$, respectively.
\end{theorem}

Note that the estimate (\ref{eq:spectral-est}) is sharp, that is,
there is an operator $T$, for which we have the equality in
(\ref{eq:spectral-est}) (see Subsection~\ref{sec:spectrum}).
\subsection{Structure of the set of periodic points}
In the following we will consider Jordan curves or arcs of two types. 
We say that an arc $\gamma$ is closed (open) if it is
homeomorphic to the segment $[0,1]$ (to the interval $(0,1)$). Let
$\Gamma$ be an oriented Jordan smooth curve. A homeomorphism
$\alpha:\Gamma\to\Gamma$ is called a shift function (shift). Put
$\alpha_0(t)=t$ and $\alpha_n(t)=\alpha[\alpha_{n-1}(t)]$ for
$n\in\bZ$ and $t\in\Gamma$.

A point $\tau\in\Gamma$ is called a periodic point 
of the multiplicity $m\ge 1$
for the shift $\alpha$, if $\alpha_m(\tau)=\tau$ and
(in the case $m>1$) $\alpha_j(\tau)\ne \tau$ for every $j=1,2,
\dots, m-1$. A periodic point of the multiplicity one is called a
fixed point. Let $\Lambda$ be the set of all periodic points. A
classification of shifts with respect to the set $\Lambda$ is
given in \cite[Ch.~1, Section~3]{KrLit}. In the following we
assume that the set $\Lambda$ is nonempty. It is known that if
$\alpha$ preserves the orientation on $\Gamma$ then all periodic
points of $\alpha$ have the same multiplicity \cite[Theorem~1.3.1]{KrLit}. 
If $\alpha$ changes the orientation on $\Gamma$ then
$\alpha$ has two fixed points on $\Gamma$, and all other periodic
points of $\alpha$ (if they exist) have the multiplicity two
\cite[Theorem 1.3.2]{KrLit}. Clearly, the set $\Lambda$ is closed.
Let $m$ be the multiplicity of periodic points of the shift
$\alpha$ if $\alpha$ preserves the orientation on $\Gamma$, and put
$m:=2$ otherwise. Then $\Lambda$ is the set of all fixed points of
the shift $\alpha_m$. Put
\[
\Phi:=\overline{\{t\in\Gamma:\alpha_m(t)\ne t\}},
\]
and decompose the contour $\Gamma$ into disjoint sets
\[
\Gamma
=
(\Gamma\setminus\Phi)\cup\Phi
=
(\Gamma\setminus\Phi)\cup(\Phi\setminus\Lambda)\cup Y
\]
where the set $Y :=\partial\Lambda=\Lambda\cap\Phi$ is the boundary of $\Lambda$.
Clearly,
\[
\Phi\setminus\Lambda=\Gamma\setminus\Lambda, \quad {\rm and} \quad
(\Gamma\setminus\Phi)\cup Y =\Lambda.
\]
The open set $\Phi\setminus\Lambda$ can be represented as a
countable (or finite) union of connected components (open arcs),
every of which is invariant with respect to the shift $\alpha_m$
($\alpha_m$-invariant) and does not contain its fixed points. The
open set $\Gamma\setminus\Phi$ also can be represented as a
countable (or finite) union of $\alpha_m$-invariant open arcs but
consisting of inner points of $\Lambda$.

Let $Y'$ denote the set of all limit points of the boundary $Y$ of
$\Lambda$.
\begin{lemma}$\!\!\!${\bf .}\label{le:decomp}
{\rm (a)} If the set $Y$ is finite, then there is a finite
decomposition
\begin{equation}\label{eq:decomp-1}
\Gamma=\left(\bigcup_i \overline{\omega_i}\right)\cup
\left(\bigcup_j \overline{\gamma_j}\right),
\end{equation}
where $\omega_i\subset\Gamma\setminus\Phi,\;
\gamma_j\subset\Phi\setminus\Lambda$ are pairwise disjoint,
$\alpha_m$-invariant open arcs with endpoints in $Y$.

{\rm (b)} Let $f:\Gamma\to\bR$ be a continuous function. If the
set $Y$ is infinite and $f(\tau)>0$ for all $\tau\in Y'$, then
there is a finite decomposition
\begin{equation}\label{eq:decomp-2}
\Gamma=\left(\bigcup_i \overline{\omega_i}\right)
\cup\left(\bigcup_j \overline{\gamma_j}\right)
\cup\left(\bigcup_k \overline{v_k}\right),
\end{equation}
where $\omega_i\subset\Gamma\setminus\Phi,
\;\gamma_j\subset\Phi\setminus\Lambda,\; v_k \subset \Gamma$ are
pairwise disjoint,  $\alpha_m$-invariant open arcs with endpoints
in $Y,\; \overline{v_k}\cap Y'\ne \emptyset$, and
$f(t)>0$ for all $t\in\overline{v_k}$ and all arcs $v_k$.
\end{lemma}
{\bf Proof.} If the set $Y$ is finite, then the set
$\Gamma\setminus Y=(\Gamma\setminus\Phi)\cup(\Phi\setminus\Lambda)$
consists of a finite set of pairwise disjoint open arcs with
endpoints in $Y$. Since the set $Y$ consists of fixed points of
$\alpha_m$, these arcs are $\alpha_m$-invariant. 
Denoting these arcs lying in $\Gamma \setminus \Phi$
and $\Phi \setminus \Lambda$  by $\omega_i$ and $\gamma_j$,
respectively, we get
\[
\Gamma\setminus\Phi=\bigcup_i\omega_i, \quad
\Phi\setminus\Lambda=\bigcup_j\gamma_j,
\]
which proves part (a).

(b) If the set $Y$ is infinite, then the set $Y'$ of its limit
points is nonempty, and vice versa. Moreover, since $Y'$ is closed
in $\Gamma$ and $\Gamma$ is compact, $Y'$ is compact as well.

Since $f$ is continuous, for every $\eps>0$ and every $\tau\in
Y'$, there is an open arc $\gamma(\tau)\ni \tau $ of the length
$|\gamma(\tau)|<\varepsilon$ and such that $f(t)>0$ for all
$t\in\overline{\gamma(\tau)}$. Without loss of generality, we can
choose endpoints of $\gamma(\tau)$ such that either the right
(left) half-neighborhood of the point $\tau$ is separated from
$Y\setminus\{\tau\}$ or the right (left) endpoint of the arc
$\gamma(\tau)$ belongs to $Y$ (in the second case the point $\tau$
is a limit point of the subset of $Y$, which lies from the right
(from the left) of the point $\tau$).

As $Y'$ is a compact set, from the open overlapping
$\{\gamma(\tau):\tau\in Y' \}$ we may extract a finite overlapping
$\{\gamma(\tau_k):k=1,2,\dots,n\}$. From each arc $\gamma(\tau_k)$
we delete half-neighborhoods of $\tau_k$ which wholly lie in
$\Lambda$ or in $(\Gamma\setminus\Lambda)\cup\{\tau_k\}$. Since
the point $\tau_k$ is not an isolated point of the set $Y$, we can
delete at most one neighborhood of $\tau_k$. Thus, from each arc
$\gamma(\tau_k)$ we obtain an arc $u(\tau_k)$, endpoints of which
belong to $Y$, and
\begin{equation}\label{eq:decomp-4}
\Big(u(\tau_k)\cap Y\Big)\cup\{\tau_k\}=\gamma(\tau_k)\cap Y.
\end{equation}
From (\ref{eq:decomp-4}) we see that the set
\[
\left(\Gamma\setminus\left(\bigcup_{k=1}^n\gamma(\tau_k)\right)\right)\cap
Y= \left(\Gamma\setminus\left(\bigcup_{k=1}^n
u(\tau_k)\right)\right)\cap Y
\]
is empty or finite.

We consecutively exclude from the set
$\{u(\tau_k):k=1,2,\dots,n\}$ the arcs $u(\tau_k)$ contained in
the union of the remaining arcs. After that we obtain a system of
arcs such that every two arcs either are not intersected,
or are intersected but each arc does not contain the other one.
Moreover, the intersection of any three such arcs is empty. Without
loss of generality, we assume that the original system
$\{u(\tau_k): k=1,2,\dots,n\}$ has this property. In that case
there are pairwise disjoint open arcs $v_k,\,k=1,2,\dots,n$, with
endpoints in $Y$, such that $ v_k\subset u(\tau_k)$ and
\begin{equation}\label{eq:decomp-5}
{\bigcup_{k=1}^n\overline{v_k}}={\bigcup_{k=1}^n\overline{u(\tau_k)}}\supset
Y'.
\end{equation}
Clearly, all arcs $v_k,\; k=1,2,\dots,n$, are
$\alpha_m$-invariant, and for every arc $v_k$,
\[
f(t)>0\quad\mbox{for all}\quad t\in\overline{v_k}, \quad\quad
\overline{v_k}\cap Y'\ne\emptyset.
\]
Further, if the set
\[
v_0:=\Gamma\setminus\left(\bigcup_{k=1}^n \overline{v_k}\right)
\]
is empty, then we have the desired decomposition of $\Gamma$.

Otherwise, the open set $v_0$ is non-empty and
$\alpha_m$-invariant. From (\ref{eq:decomp-5}) we see that
\[
v_0\cap Y= \left(\Gamma\setminus\left(\bigcup_{k=1}^n\overline
{u(\tau_k)}\right)\right)\cap( Y\setminus Y'),
\]
that is, the set $v_0\cap Y$ is finite.

Since the open set $v_0$ is not empty, the set
\[
v_0\cap\Big(\Gamma\setminus
Y\Big)=\Big(v_0\cap(\Gamma\setminus\Phi)\Big)\cup\Big(v_0\cap(\Phi\setminus\Lambda)\Big)
\]
consists of a finite set of
open arcs with endpoints in $v_0\cap Y$. Denoting these arcs by
$\omega_i$ and $\gamma_j$ if they lie, respectively, in
$\Gamma\setminus\Phi$ and $\Phi\setminus\Lambda$, we obtain
\[
v_0\cap(\Gamma\setminus\Phi)=\bigcup_i\omega_i, \quad
v_0\cap(\Phi\setminus\Lambda)=\bigcup_j\gamma_j.
\]
Clearly, these arcs are $\alpha_m$-invariant, and we get the desired decomposition
(\ref{eq:decomp-2}). \rule{2mm}{2mm}

\begin{proposition}$\!\!\!${\bf .}\label{pr:isolated}
The set $Y\setminus Y'$ of all isolated points of $Y$ is contained
in the set of all endpoints of connected components
$\gamma\subset\Phi\setminus\Lambda(=\Gamma\setminus\Lambda)$.
\end{proposition}

The proof is obvious.
\section{Two-sided invertibility: sufficient conditions}
\setcounter{equation}{0}
\subsection{Estimate for the spectral radius of weighted shift operators}
In Sections 3 and 4 we suppose that $\Gamma$ is an oriented Jordan smooth curve and $\alpha$
is an orientation preserving  diffeomorphism of $\Gamma$ onto itself which has
an arbitrary nonempty set $\Lambda$ of fixed points.
We equip $\Gamma$ with the Lebesgue length measure $|d\tau|$. In the following
we will consider all rearrangement-invariant spaces over the finite measure
space $(\Gamma,|d\tau|)$. Also we always suppose that r.-i. spaces $X(\Gamma)$
have nontrivial Boyd indices $\alpha_X,\beta_X$. In that case the operator $A=aI-bW$
is bounded in the space $X(\Gamma)$, due to the Boyd interpolation theorem \cite{b4}.

If $\Omega$ is a subset of $\Gamma$ of a positive measure, then
we will denote by $X(\Omega)$ the subspace of $X(\Gamma)$ which
consists of functions supported in $\Omega$. In view of
Lemma~\ref{le:subset}, $X(\Omega)$ is an r.-i. space with the same
Boyd and Zippin indices as the whole space $X(\Gamma)$.

Denote by $C(\Gamma)$ the set of all continuous functions on $\Gamma$.
\begin{theorem}$\!\!\!${\bf .}\label{th:spectr-Lp}
The spectral radius of the weighted shift operator $gW$ with
coefficient (weight) $g\in C(\Gamma)$ in the Lebesgue space
$L^p(\Gamma), 1< p<\infty$, is calculated by the formula
\[
r(gW;L^p(\Gamma))=\max_{\tau\in\Lambda} |g(\tau)||\alpha'(\tau)|^{-1/p}.
\]
\end{theorem}

This theorem follows from the results of \cite{KK76}. For further
generalizations of the formula for the spectral radius of the
weighted shift operator, see \cite[Ch.~1, Section~5]{Antonevich}.
\begin{theorem}$\!\!\!${\bf .}\label{th:spectral-X}
The spectral radius of the weighted shift operator $gW$ with
coefficient (weight) $g\in C(\Gamma)$ in an r.-i. space
$X(\Gamma)$ with nontrivial Boyd indices $\alpha_X,\beta_X$
satisfies the estimate
\[
r(gW;X(\Gamma))\le\max_{\tau\in\Lambda} \Big(|g(\tau)|\max
\Big\{|\alpha'(\tau)|^{-\alpha_X},|\alpha'(\tau)|^{-\beta_X}\Big\}\Big).
\]
\end{theorem}
{\bf Proof.}
By Theorem~\ref{th:spectr-Lp}, the spectral radii of the operator $gW$
in the Lebesgue spaces $L^{1/\alpha_X}(\Gamma)$ and $L^{1/\beta_X}(\Gamma)$,
are calculated by
\begin{eqnarray*}
r(gW;L^{1/\alpha_X}(\Gamma)) =
\max_{\tau\in\Lambda}|g(\tau)||\alpha'(\tau)|^{-\alpha_X},
\quad
r(gW;L^{1/\beta_X}(\Gamma)) =
\max_{\tau\in\Lambda}|g(\tau)||\alpha'(\tau)|^{-\beta_X},
\end{eqnarray*}
respectively. Theorem~\ref{th:spectral-est} and the latter equalities  give
\begin{eqnarray*}
r(gW;X(\Gamma))
&\le&
\max\Big\{r(gW;L^{1/\alpha_X}(\Gamma)),r(gW;L^{1/\beta_X}(\Gamma))\Big\}
\\
&=&
\max\left\{
\max_{\tau\in\Lambda}|g(\tau)||\alpha'(\tau)|^{-\alpha_X},
\max_{\tau\in\Lambda}|g(\tau)||\alpha'(\tau)|^{-\beta_X}
\right\}
\\
&=&
\max_{\tau\in\Lambda}
 \Big(|g(\tau)|\max
\Big\{|\alpha'(\tau)|^{-\alpha_X},|\alpha'(\tau)|^{-\beta_X}\Big\}
\Big). \quad\rule{2mm}{2mm}
\end{eqnarray*}
\subsection{Sufficient conditions for the two-sided invertibility}
For the functional operator (\ref{eq:operator})
define the two continuous functions $\eta_i:\Gamma\to\bR, i\in\{0,1\}$, by the formulas
\begin{eqnarray}
\label{eq:eta0}
\eta_0(t)&:=&|a(t)|-|b(t)|\min\Big\{|\alpha'(t)|^{-\alpha_X},|\alpha'(t)|^{-\beta_X}\Big\},\\
\label{eq:eta1}
\eta_1(t)&:=&|a(t)|-|b(t)|\max\Big\{|\alpha'(t)|^{-\alpha_X},|\alpha'(t)|^{-\beta_X}\Big\}.
\end{eqnarray}
\begin{lemma}$\!\!\!${\bf .}\label{le:invert}
Let $\gamma$ be a closed $\alpha$-invariant arc of $\Gamma$.

\noindent
{\rm (a)} If $\eta_0(t)<0$ for every $t\in\gamma$,
then the operator $A$ is two-sided invertible in 
$X(\gamma)$, and
\begin{equation}\label{eq:invert-1}
A^{-1}=-W^{-1}\sum_{n=0}^\infty(b^{-1}aW^{-1})^nb^{-1}I.
\end{equation}

\noindent
{\rm (b)} If $\eta_1(t)>0$ for every $t\in\gamma$,
then the operator $A$ is two-sided invertible in 
$X(\gamma)$, and
\[
A^{-1}=\sum_{n=0}^\infty(a^{-1}bW)^na^{-1}I.
\]
\end{lemma}
{\bf Proof.} We consider the more difficult part (a). From the
definition of $\eta_0$ we see that if $\eta_0(t)<0$ for all
$t\in\gamma$, then
\begin{equation}\label{eq:invert-2}
0\le |a(t)|<|b(t)|\min\Big\{|\alpha'(t)|^{-\alpha_X},|\alpha'(t)|^{-\beta_X}\Big\},
\quad t\in\gamma.
\end{equation}
Hence, $b$ is invertible in $C(\gamma)$. In that case
\begin{equation}\label{eq:invert-3}
A=-b(I-b^{-1}aW^{-1})W.
\end{equation}
Let us show the invertibility of the factor $I-b^{-1}aW$. 
For fixed points of $\alpha$ we get
$\alpha'_{-1}(\tau)=1/\alpha'[\alpha_{-1}(\tau)]=1/\alpha'(\tau).$
Then, taking into account that $W_\alpha^{-1}=W_{\alpha_{-1}}$ and the arc $\gamma$ is
$\alpha$-invariant, we infer from
Theorem~\ref{th:spectral-X} that
\begin{eqnarray}
r(b^{-1}aW^{-1};X(\gamma))
&\le&
\max_{\tau\in\Lambda\cap\gamma}
\left(
|b(\tau)|^{-1}|a(\tau)|\max\Big\{|\alpha_{-1}'(\tau)|^{-\alpha_X},|\alpha_{-1}'(\tau)|^{-\beta_X}
\Big\}
\right)
\nonumber
\\
& =&
\max_{\tau\in\Lambda\cap\gamma}
\left(
|a(\tau)|\Big(|b(\tau)|\min\Big\{|\alpha'(\tau)|^{-\alpha_X},|\alpha'(\tau)|^{-\beta_X}
\Big\}\Big)^{-1}
\right).
\label{eq:invert-4}
\end{eqnarray}
From (\ref{eq:invert-2}) and (\ref{eq:invert-4}) we see that
$r(b^{-1}aW^{-1};X(\gamma))<1$. Hence, the operator
$I-b^{-1}aW^{-1}$ is invertible in $X(\gamma)$, and
\begin{equation}\label{eq:invert-5}
(I-b^{-1}aW^{-1})^{-1}=\sum_{n=0}^\infty(b^{-1}aW^{-1})^n.
\end{equation}
From (\ref{eq:invert-3}), (\ref{eq:invert-5}) and the
invertibility of operators $bI$ and $W$ we conclude that $A$ is
invertible and (\ref{eq:invert-1}) holds. The assertion (b) can be
proved analogously. \rule{2mm}{2mm}
\section{One-sided invertibility: case of fixed points}
\setcounter{equation}{0}
\subsection{Case of only two fixed points}
Let $X(\Gamma)$ be a reflexive r.-i. space of fundamental type
with nontrivial Boyd indices, i.e., Boyd and Zippin indices
satisfy the inequalities
\[
0<\alpha_X=p_X\le q_X=\beta_X<1.
\]

We say that an $\alpha$-invariant set $\Omega\subset\Gamma$
satisfies, respectively, the condition
\begin{eqnarray*}
R(\Omega)& {\rm if}&\mbox{for every $t\in\Omega$ there exists a
$k_0=k_0(t)\in\bZ$ such that }
\\
&& a[\alpha_k(t)]\ne 0 \mbox{ for  } k<k_0 \quad \mbox{ and }
\quad b[\alpha_k(t)]\ne 0 \mbox{ for } k\ge k_0,\\ L(\Omega)&{\rm
if}&\mbox{for every $t\in\Omega$ there exists a $k_0=k_0(t)\in\bZ$
such that }
\\
&&
a[\alpha_k(t)]\ne 0
\mbox{ for  }
k>k_0
\quad
\mbox{ and }
\quad
b[\alpha_k(t)]\ne 0
\mbox{ for }
k< k_0.
\end{eqnarray*}

Denote by $GC(\Omega)$ the set of all functions invertible in
$C(\Omega)$, that is, the set of all functions $a\in C(\Omega)$
such that
\[
\inf_{t\in\Omega}|a(t)|>0.
\]

In the case of only two fixed points we obtain from  \cite[Theorem~6.8]{AKarl00}
and Lemma~\ref{le:subset} the following result.
\begin{theorem}$\!\!\!${\bf .}\label{th:two-fixed}
Let $\gamma\subset\Gamma$ be an open arc with endpoints $\tau_\pm$
and $\Lambda\cap\overline{\gamma}=\{\tau_-,\tau_+\}$.
The operator $A$ is right (left) invertible in the subspace
$X(\gamma)$ if and only if one of the three conditions holds:
\begin{equation}
\label{eq:two-fixed1}
\eta_1(\tau_-)>0,\quad\eta_1(\tau_+)>0,\quad a\in GC(\overline{\gamma});
\end{equation}
\begin{equation}
\label{eq:two-fixed2}
\eta_0(\tau_-)<0,\quad\eta_0(\tau_+)<0,\quad b\in GC(\overline{\gamma});
\end{equation}
\[
\eta_0(\tau_+)<0<\eta_1(\tau_-)\quad and\quad R(\gamma)\mbox{ is
fulfilled}
\]
(respectively, {\rm (\ref{eq:two-fixed1}), (\ref{eq:two-fixed2})}, or
\begin{equation}\label{eq:two-fixed4}
\eta_0(\tau_-)<0<\eta_1(\tau_+)\quad and \quad L(\gamma)\mbox{ is
fulfilled }  \mbox{)}.
\end{equation}
\end{theorem}

The restrictions $\alpha_X=p_X$ and $\beta_X=q_X$ are essential
for the proof of this theorem given in \cite{AKarl00}.
From Theorem~\ref{th:two-fixed} and the inequality
$\eta_1(t)\le\eta_0(t),t\in\Gamma$, it follows
\begin{corollary}$\!\!\!${\bf .}\label{co:invert-one-sided}
If the operator $A$ is one-sided invertible in the subspace
$X(\gamma)$, then
\[
\eta_0(\tau)\eta_1(\tau)>0,\quad\tau\in\{\tau_-,\tau_+\}.
\]
\end{corollary}
\begin{corollary}$\!\!\!${\bf .}\label{le:onesided-spectrum}
The intersection of the left and right spectra of the shift
operator $W$ in the subspace $X(\gamma)$, that is, the set
\[
\sigma_0 (W):= \Big\{ \lambda\in\bC\ :\ \lambda I-W\quad\mbox{is
not right and not left invertible in } X(\gamma) \Big\},
\]
is given by the formula
\[
\sigma_0 (W)=\bigcup_{\tau\in\{\tau_-,\tau_+\}}
\Big\{\lambda\in\bC\ :\ \delta(\tau)\le |\lambda|\le \Delta(\tau)
\Big\},
\]
where
$
\delta(\tau)
:=
\min\Big\{|\alpha'(\tau)|^{-\alpha_X},|\alpha'(\tau)|^{-\beta_X}\Big\},
\quad
\Delta(\tau)
:=
\max\Big\{|\alpha'(\tau)|^{-\alpha_X},|\alpha'(\tau)|^{-\beta_X}\Big\}.
$
\end{corollary}
{\bf Proof.}
By Corollary~\ref{co:invert-one-sided}, if the operator
$\lambda I-W$ is one-sided invertible in $X(\gamma)$, then
\[
\eta_0(\tau)\eta_1(\tau)=(|\lambda|-\delta(\tau))(|\lambda|-\Delta(\tau))
>0
\quad\mbox{for}\quad\tau\in\{\tau_-,\tau_+\},
\]
or, equivalently,
\begin{equation}\label{eq:onesided-spectrum-1}
|\lambda|\not\in[\delta(\tau),\Delta(\tau)]
\quad\mbox{for}\quad\tau\in\{\tau_-,\tau_+\}.
\end{equation}
On the other hand, if (\ref{eq:onesided-spectrum-1}) does not hold, then
one of the following four conditions is satisfied:
\begin{tabbing}
xxx\=xxxxxxxxxxxxxxxxxxxxxx\=xxxxx\=xxxxxxxxxxxxxxxx\kill
1) \> $|\lambda|<\min\{\delta(\tau_-),\delta(\tau_+)\}$
\> $\Longleftrightarrow$
\> $\eta_0(\tau_-)<0$ and $\eta_0(\tau_+)<0$,
\\
2)
\> $|\lambda|>\max\{\Delta(\tau_-),\Delta(\tau_+)\}$
\> $\Longleftrightarrow$
\> $\eta_1(\tau_-)>0$ and $\eta_1(\tau_+)>0$,
\\
3)
\> $\Delta(\tau_-)<|\lambda|<\delta(\tau_+)$
\> $\Longleftrightarrow$
\> $\eta_0(\tau_+)<0<\eta_1(\tau_-)$,
\\
4)
\> $\Delta(\tau_+)<|\lambda|<\delta(\tau_-)$
\> $\Longleftrightarrow$
\> $\eta_0(\tau_-)<0<\eta_1(\tau_+)$.
\end{tabbing}
Since $\lambda\ne 0$ in the cases 2)--4), we infer from Theorem~\ref{th:two-fixed}
that the operator $\lambda I-W$ is two-sided invertible in $X(\gamma)$ in cases
1)--2), right invertible in case 3), and left invertible in case 4).

Thus, the operator $\lambda I-W$ is one-sided invertible in the subspace $X(\gamma)$ if
and only if (\ref{eq:onesided-spectrum-1}) holds. Hence, $\lambda I-W$ is neither
left invertible, nor right invertible in the subspace  $X(\gamma)$ if and only 
if $\lambda\in\sigma_0(W)$.
\rule{2mm}{2mm}

Clearly, if $\alpha_X<\beta_X$, then $\sigma_0$ is the union of
two annuli which degenerate in two circles whenever
$\alpha_X=\beta_X$. Thus, in the case of arbitrary r.-i. spaces,
in contrast to Lebesgue spaces, the following new quality appears:
the intersection of the left and right spectra of the shift
operator $W$ becomes massive (i.e., admits non-zero plain measure)
in general.
\subsection{Case of arbitrary nonempty set of fixed points}
In this subsection we prove a criterion for one-sided
invertibility of the operator $A$ when the shift $\alpha$ has  an
arbitrary nonempty set of fixed points. In the proof we
essentially use a decomposition of $\Gamma$ into a union of arcs
of three types (see Lemma~\ref{le:decomp}). Since these arcs are
invariant with respect to the shift $\alpha$, we can check
one-sided invertibility of $A$ in each subspace of functions
supported in these arcs.
\begin{theorem}$\!\!\!${\bf .}\label{th:fixed}
The operator $A$ is right
(left) invertible in the space $X(\Gamma)$ if and only if the following
three conditions simultaneously hold:

\noindent
{\rm (i)} the operator $A$ is right (left) invertible in 
$X(\gamma)$ for every connected component
$\gamma\subset\Phi\setminus\Lambda$; 

\noindent
{\rm (ii)} $a(t)\ne b(t)$ for every $t\in\Gamma\setminus\Phi$;

\noindent
{\rm (iii)} if the set $Y$ is infinite, then
$\eta_0(\tau)\eta_1(\tau)>0$ for every $\tau\in Y'$.
\end{theorem}
{\bf Proof.} {\it Necessity.} Suppose the operator $A$ is right
(left) invertible in $X(\Gamma)$. Since the set
$\Gamma\setminus\Phi$ and every connected component
$\gamma\subset\Phi\setminus\Lambda$ are $\alpha$-invariant, the
operator $A$ is bounded and right (left) invertible in every
subspace $X(\gamma),\gamma\subset\Phi\setminus\Lambda$, and in the
subspace $X(\Gamma\setminus\Phi)$. So, we get (i).

Since $\Gamma\setminus\Phi=\{t\in\Gamma\setminus Y:\alpha(t)=t\}$,
we have $A=(a-b)I$ in the subspace $X(\Gamma\setminus\Phi)$. In
view of the one-sided invertibility of $(a-b)I$, we obtain (ii).
If the set $Y$ is finite, the necessity is proved.

Now we prove (iii) when the set $Y$ is infinite. Assume the
contrary: $\eta_0(\tau_0)\eta_1(\tau_0)\le 0$ for some point
$\tau_0\in Y'$. In view of the stability of one-sided
invertibility of an operator under small (in the operator norm)
perturbations (see, e.g., \cite[Ch.~2, Theorem~5.4]{gk}), there
are $\widetilde{a},\widetilde{b}\in C(\Gamma)$ such that the
operator $\widetilde{a}I-\widetilde{b}W$ is right (left)
invertible in $X(\Gamma)$, and at some point $\tau\in Y\setminus
Y'$, which is sufficiently close to $\tau_0$, the following
inequality holds:
$\widetilde{\eta}_0(\tau)\widetilde{\eta}_1(\tau)\le 0$. Here
quantities $\widetilde{\eta}_0$ and $\widetilde{\eta}_1$ are
defined for the operator $\widetilde{a}I-\widetilde{b}W$ by
formulas (\ref{eq:eta0}) and (\ref{eq:eta1}).
Proposition~\ref{pr:isolated} tells us that $\tau$ is an endpoint of
some connected component $\gamma\subset\Phi\setminus\Lambda$. From
Corollary~\ref{co:invert-one-sided} we obtain that the operator
$\widetilde{a}I-\widetilde{b}W$ is not one-sided invertible in
$X(\gamma)$. This contradicts to condition (i) just proved. Hence,
condition (iii) is fulfilled. Necessity is proved.

{\it Sufficiency.} Since by (i) the operator $A$ is right (left)
invertible in $X(\gamma)$ for every connected component
$\gamma\subset \Phi\setminus\Lambda$, it follows from Corollary
\ref{co:invert-one-sided} that $\eta_0(\tau)\eta_1(\tau)>0$ for
each endpoint $\tau$ of every such $\gamma$. Hence, by Proposition
\ref{pr:isolated},
\[
\eta_0(\tau)\eta_1(\tau)>0 \quad {\rm for} \; {\rm every} \quad
\tau\in Y\setminus Y',
\]
which together with (iii) gives
\begin{equation}\label{f4.4'}
\eta_0(\tau)\eta_1(\tau)>0 \quad {\rm for} \; {\rm every} \quad
\tau\in Y.
\end{equation}

Since $\alpha'(\tau)=1$ in every endpoint $\tau$ of every
connected component $\omega\subset \Gamma\setminus \Phi$, we infer
from (\ref{eq:eta0})-(\ref{eq:eta1}) and (\ref{f4.4'}) that
$a(\tau)\ne b(\tau)$ for all those points $\tau$. Then in view of
(ii), for every connected component $\omega\subset\Gamma\setminus\Phi$,
\begin{equation}\label{f4.4''}
  a(t) \ne b(t) \quad {\rm for} \; {\rm every} \quad t \in
  \overline {\omega}.
\end{equation}

If the set $Y$ is finite, then condition (iii)
disappears. Moreover, by Lemma~\ref{le:decomp}(a) there is a
finite decomposition
\[
\Gamma=
\left(\bigcup_i\overline{\omega_i}\right)\cup
\left(\bigcup_j\overline{\gamma_j}\right),
\]
where $\omega_i\subset \Gamma\setminus\Phi$ and
$\gamma_j\subset\Phi\setminus \Lambda$ are pairwise disjoint,
$\alpha$-invariant open arcs with endpoints in $Y$. Decompose the
space $X(\Gamma)$ into the (finite) direct sum of its subspaces:
\begin{equation}\label{eq:fixed-1}
X(\Gamma)=\left(\bigoplus_i X(\omega_i)\right)\oplus
\left(\bigoplus_j X(\gamma_j)\right).
\end{equation}
Clearly, the operator $A$ acts in each subspace of the
decomposition (\ref{eq:fixed-1}). Since
$\omega_i\subset\Gamma\setminus\Phi=\{t\in\Gamma\setminus Y
:\alpha(t)=t\}$, from (\ref{f4.4''}) we see that the operator
$A=(a-b)I$ is invertible in each subspace $X(\omega_i)$. Hence,
taking into account condition (i), we obtain that the operator $A$
is right (left) invertible in the space $X(\Gamma)$. Thus, in the
case of finite $Y$ sufficiency is proved.

Now we consider the case of infinite set $Y$. Since the functions
$\eta_0$ and $\eta_1$ are continuous on $\Gamma$, from
Lemma~\ref{le:decomp}(b) and condition (iii) it follows that there
is a finite decomposition
\[
\Gamma=\left(\bigcup_i\overline{\omega_i}\right) \cup
\left(\bigcup_j\overline{\gamma_j}\right) \cup
\left(\bigcup_k\overline{v_k}\right),
\]
where $\omega_i\subset \Gamma\setminus \Phi,\;\gamma_j\subset
\Phi\setminus \Lambda$ and $v_k\subset \Gamma$ are pairwise
disjoint, $\alpha$-invariant open arcs with endpoints in $Y$.
Moreover, $\overline{v_k}\cap Y' \ne \emptyset$ and
\begin{equation}\label{eq:fixed-2}
\eta_0(\tau)\eta_1(\tau)>0\quad\mbox{for all}\quad t\in\overline{v_k}
\end{equation}
and all arcs $v_k$. Decompose the space $X(\Gamma)$ into the (finite)
direct sum of its subspaces:
\begin{equation}\label{eq:fixed-3}
X(\Gamma)=
\left(\bigoplus_i X(\omega_i)\right)
\oplus\left(\bigoplus_j X(\gamma_j)\right)
\oplus\left(\bigoplus_k X(v_k)\right).
\end{equation}
Clearly, the operator $A$ acts in each subspace of the
decomposition (\ref{eq:fixed-3}). From (\ref{eq:fixed-2}) we see
that either $\eta_0(\tau)\ge\eta_1(\tau)>0$, or
$\eta_1(\tau)\le\eta_0(\tau)<0$ for every $\tau\in\overline{v_k}$.
Hence, by Lemma~\ref{le:invert}, the operator $A$ is invertible in
each subspace $X(v_k)$. As in the previous case, we see that the
operator $A$ is invertible in each  subspace $X(\omega_i)$ and is
right (left) invertible in each subspace $X(\gamma_j)$. Thus,
taking into account (\ref{eq:fixed-3}), we see that the operator
$A$ is right (left) invertible in the space $X(\Gamma)$.
\rule{2mm}{2mm}
\subsection{Another form of the criterion for one-sided invertibility}
In this subsection we reformulate results of the previous
subsection in terms of a function, which controls the two- and
one-sided invertibility of the operator $A$. We start with the
following important property of fixed points.
\begin{lemma}\label{le:atr-repel}
{\rm (see \cite[Ch.~1, Lemma~2]{KrLit}).} Let
$\tau_1,\tau_2\in\Gamma$ be fixed points of an orientation
preserving shift $\alpha$. If the arc $(\tau_1,\tau_2)$ does not
contain fixed points of the shift $\alpha$, then for each point
$t\in(\tau_1,\tau_2)$ the iterative sequence $\alpha_n(t)\
(\alpha_{-n}(t) )$ converges to a fixed point, either to the point
$\tau_1$ (respectively, $\tau_2$) or to the point $\tau_2$
(respectively, $\tau_1$) independently of $t$.
\end{lemma}

In accordance with this property we define the following
functions:
\begin{equation}\label{eq:etapm-def}
\eta_i^\pm(t):=\lim_{n\to\pm\infty}\eta_i[\alpha_n(t)],
\quad t\in\Gamma, \quad i\in\{0,1\}.
\end{equation}
In view of Lemma~\ref{le:atr-repel} and the continuity of
$a,b,\alpha'$, the limits in (\ref{eq:etapm-def}) exist for every
$t\in\Gamma$. Hence the functions $\eta_i^\pm$ are well-defined.
From the definition of fixed points we see that
\begin{equation}\label{eq:etapm}
\eta_i(\tau)=\eta_i^+(\tau)=\eta_i^-(\tau),
\quad\tau\in\Lambda, \quad i\in\{0,1\}.
\end{equation}

For the functional operator (\ref{eq:operator}), define five sets:
$\Gamma_1:=\Gamma\setminus\Phi$,
\begin{eqnarray}\label{eq:Gamma-def1}
\Gamma_2 :=\Big\{t\in\Phi:\eta_1^-(t)>0,\ \eta_1^+(t)>0\Big\},
&&
\Gamma_3 :=\Big\{t\in\Phi:\eta_0^-(t)<0,\ \eta_0^+(t)<0\Big\},
\\
\label{eq:Gamma-def2} \Gamma_4
:=\Big\{t\in\Phi:\eta_0^+(t)<0<\eta_1^-(t)\Big\},\quad\: &&
\Gamma_5 :=\Big\{t\in\Phi:\eta_0^-(t)<0<\eta_1^+(t)\Big\}.
\end{eqnarray}
Clearly, these sets are pairwise disjoint and $\alpha$-invariant.
Moreover, we infer from (\ref{eq:etapm}), the inequality
$\eta_1(\tau)\le\eta_0(\tau)$
 and the definitions of $\Gamma_4$ and
$\Gamma_5$, that
\begin{equation}\label{eq:Lambda-not-in45}
\Gamma_4\cap\Lambda=\Gamma_5\cap\Lambda=\emptyset.
\end{equation}
Indeed, if, for example, $\tau\in \Gamma_4\cap \Lambda$, then
$
\eta_1^+(\tau)\le\eta_0^+(\tau)<0<\eta_1^-(\tau)\le\eta_0^-(\tau),
$
which is impossible due to (\ref{eq:etapm}).
Define the function
\[
\sigma_A(t):=\left\{
\begin{array}{cl}
a(t)-b(t), & t\in\Gamma_1,\\
a(t),      & t\in\Gamma_2,\\
-b(t),     & t\in\Gamma_3,\\
0,         & t\in\Gamma\setminus(\Gamma_1\cup\Gamma_2\cup\Gamma_3).
\end{array}
\right.
\]
Now reformulate Theorem~\ref{th:fixed} in terms of the function
$\sigma_A$ and sets $\Gamma_j$, $j\in \{1,2, \ldots ,5 \}$.
\begin{theorem}$\!\!\!${\bf .}\label{th:fixed2}
The operator $A$ is right (left) invertible in the space $X(\Gamma)$ if and only if
\begin{equation}\label{eq:fixed2-right}
\sigma_A(t)\ne 0\mbox{ for all }t\in\Gamma\setminus\Gamma_4
\quad\mbox{ and }\quad R(\Gamma_4)\mbox{ holds}
\end{equation}
(respectively,
\begin{equation}\label{eq:fixed2-left}
\sigma_A(t)\ne 0\mbox{ for all }t\in\Gamma\setminus\Gamma_5
\quad\mbox{ and }\quad L(\Gamma_5)\mbox{ holds )}.
\end{equation}
\end{theorem}
{\bf Proof.} {\it Necessity}. Let the operator $A$ be one-sided
invertible in the space $X(\Gamma)$. Then conditions (i)--(iii) of
Theorem~\ref{th:fixed} are satisfied. The proof of sufficiency in
that theorem shows that conditions (i) and (iii) imply
(\ref{f4.4'}). Hence, for $\tau\in Y$, either
$\eta_0(\tau)\ge\eta_1(\tau)>0$, or
$\eta_1(\tau)\le\eta_0(\tau)<0$, which gives $Y\subset \Gamma_2
\cup \Gamma_3$. Taking into account (\ref{eq:eta1}) and
(\ref{eq:etapm}), we conclude that in the first case
$\tau\in\Gamma_2\cap Y$ and then
\[
|a(\tau)|>|b(\tau)|
\max\Big\{|\alpha'(\tau)|^{-\alpha_X},|\alpha'(\tau)|^{-\beta_X}\Big\}\ge 0.
\]
Hence,
\begin{equation}\label{eq:fixed2-3}
\sigma_A(\tau)\ne 0\quad\mbox{for all}\quad\tau\in\Gamma_2\cap Y.
\end{equation}
Analogously, in the second case $\tau\in \Gamma_3 \cap Y$ and
\begin{equation}\label{eq:fixed2-4}
\sigma_A(\tau)\ne 0\quad\mbox{for all}\quad\tau\in\Gamma_3\cap Y.
\end{equation}
On the other hand, condition (ii) of Theorem~\ref{th:fixed} implies
\begin{equation}\label{eq:fixed2-5}
\sigma_A(t)\ne 0\quad\mbox{for all}\quad t\in\Gamma_1.
\end{equation}

Consider the case of the right invertibility. By
Theorem~\ref{th:fixed}(i), the operator $A$ is right invertible in
the subspace $X(\gamma)$ for every connected component
$\gamma\subset\Phi\setminus\Lambda$. Hence,
Theorem~\ref{th:two-fixed} yields that one of the three conditions
holds:
\begin{eqnarray}
&&
\gamma\subset(\Phi\setminus\Lambda)\cap\Gamma_2
\quad\mbox{and}\quad\sigma_A(t)\ne 0,\quad t\in\gamma;
\label{eq:fixed2-6}
\\
&&
\gamma\subset(\Phi\setminus\Lambda)\cap\Gamma_3
\quad\mbox{and}\quad\sigma_A(t)\ne 0,\quad t\in\gamma;
\label{eq:fixed2-7}
\\
&& \gamma\subset\Gamma_4\quad\mbox{and}\quad R(\gamma)\quad
\mbox{is fulfilled}. \label{eq:fixed2-8}
\end{eqnarray}
Clearly, $\gamma\cap\Gamma_5=\emptyset$. Thus,
\[
\Gamma=\Gamma_1\cup\Gamma_2\cup\Gamma_3\cup\Gamma_4, \quad
Y\subset\Gamma_2\cap\Gamma_3,\quad\Gamma_5=\emptyset,
\]
and (\ref{eq:fixed2-3})--(\ref{eq:fixed2-7}) imply that $\sigma_A(t)\ne 0$
for all $t\in\Gamma\setminus\Gamma_4=\Gamma_1\cup\Gamma_2\cup\Gamma_3$.
Moreover, from (\ref{eq:fixed2-8}) we see that $R(\Gamma_4)$ holds.

Analogously, in the case of the left invertibility we have
\[
\Gamma=\Gamma_1\cup\Gamma_2\cup\Gamma_3\cup\Gamma_5, \quad
Y\subset\Gamma_2\cap\Gamma_3,\quad\Gamma_4=\emptyset,
\]
$\sigma_A(t)\ne 0$ for all $t\in\Gamma\setminus\Gamma_5=
\Gamma_1\cup\Gamma_2\cup\Gamma_3$ and $L(\Gamma_5)$ holds.
Necessity is proved.

{\it Sufficiency.} Suppose (\ref{eq:fixed2-left}) is fulfilled.
Since $\sigma_A(t)\ne 0$ for all $t\in\Gamma\setminus\Gamma_5$, we
see that
$\Gamma\setminus\Gamma_5=\Gamma_1\cup\Gamma_2\cup\Gamma_3$. Then
from (\ref{eq:Lambda-not-in45}) and the definition of $\Gamma_1$
we infer that $Y\subset\Gamma_2\cup\Gamma_3$. If $\tau\in
Y\cap\Gamma_2$, then $\eta_0(\tau)\ge\eta_1(\tau)>0$. Analogously,
$\eta_1(\tau)\le\eta_0(\tau)<0$ whenever $\tau\in Y\cap\Gamma_3$.
Hence, $\eta_0(\tau)\eta_1(\tau)>0$ for all $\tau\in Y$. In
particular, if the set $Y$ is infinite, then taking into account
the embedding $Y'\subset Y$, we get condition (iii) of
Theorem~\ref{th:fixed}. The condition $\sigma_A(t)\ne 0$ for all
$t\in\Gamma_1$, implies condition (ii) of Theorem~\ref{th:fixed}.

From Lemma~\ref{le:atr-repel} it follows that the set
$(\Phi\setminus\Lambda)\cap\Gamma_2,
\Big((\Phi\setminus\Lambda)\cap\Gamma_3,(\Phi\setminus\Lambda)\cap\Gamma_5\Big)$
consists of at most countable union of connected components (open
arcs), every of which wholly lies in
$(\Phi\setminus\Lambda)\cap\Gamma_2$ (respectively, in
$(\Phi\setminus\Lambda)\cap\Gamma_3,(\Phi\setminus\Lambda)\cap\Gamma_5$).
Moreover, each such arc $\gamma$ does not contain fixed points of
the shift $\alpha$, and the endpoints of this arc lie in $Y$.
Hence, for every such arc we can apply Theorem~\ref{th:two-fixed}
(case (\ref{eq:two-fixed1}), (\ref{eq:two-fixed2}),
(\ref{eq:two-fixed4}), respectively). By
Theorem~\ref{th:two-fixed}, the operator $A$ is left invertible in
the subspace $X(\gamma)$ for every connected component $\gamma$ of
the set
\[
\Phi\setminus\Lambda=
\Big((\Phi\setminus\Lambda)\cap\Gamma_2\Big)\cup
\Big((\Phi\setminus\Lambda)\cap\Gamma_3\Big)\cup
\Big((\Phi\setminus\Lambda)\cap\Gamma_5\Big).
\]
Hence, condition (i) of Theorem~\ref{th:fixed} is satisfied along
with conditions (ii) and (iii). Therefore, we obtain the required
result from Theorem~\ref{th:fixed}.

Analogously, one can derive that (\ref{eq:fixed2-right}) implies the right
invertibility of $A$.
\rule{2mm}{2mm}
\section{One-sided invertibility: general case}
\setcounter{equation}{0}
\subsection{Main result}
In this section we generalize results of the previous section to
the case of shifts having arbitrary nonempty sets of periodic
points. Let $\Gamma$ be an oriented Jordan smooth curve. Suppose
$\alpha$ is a diffeomorphism of $\Gamma$ onto itself which
preserves or changes the orientation on $\Gamma$ and has an arbitrary
nonempty set $\Lambda$ of periodic points of multiplicity $m\ge 1$
if $\alpha$ preserves the orientation on $\Gamma$, and the multiplicity
$m=2$ (for all points of $\Lambda$ except for two fixed points of
$\alpha$) if $\alpha$ changes the orientation on $\Gamma$. Suppose
$X(\Gamma)$ is a reflexive r.-i. space of fundamental type with
nontrivial Boyd indices.

For a continuous function $f:\Gamma\to\bC$ and all $k\in\bN$, we
introduce the functions
\[
f_k(t):=\prod_{i=0}^{k-1}f[\alpha_i(t)],\quad t\in \Gamma .
\]
This notation is consistent with the formula for the derivative of the shift $\alpha_k$
if $f=\alpha'$.
\begin{proposition}$\!\!\!${\bf .}\label{pr:lim}
Suppose $\alpha$ is a shift having periodic points of multiplicity
$m\ge 1$ if $\alpha$ preserves the orientation on $\Gamma$, and let
$m=2$ if $\alpha$ changes the orientation on $\Gamma$. If
$f:\Gamma\to\bC$ is a continuous function, then for every
$t\in\Gamma$ and every $k\in\bZ$,
\begin{equation}\label{eq:lim-1}
\lim_{n\to\pm\infty}
 f_m[\alpha_{mn+k}(t)]
=
\lim_{n\to\pm\infty}
f_m[\alpha_{mn}(t)].
\end{equation}
\end{proposition}
{\bf Proof.} Clearly, it is sufficient to consider the case $1\le
k\le m-1$. Then
\begin{equation}\label{eq:lim-2}
\left( \prod_{i=0}^{k-1} f[\alpha_{mn+i}(t)] \right)
 f_m[\alpha_{mn+k}(t)]=
\left( \prod_{i=0}^{k-1} f[\alpha_{m(n+1)+i}(t)] \right)
f_m[\alpha_{mn}(t)].
\end{equation}
Since the shift $\alpha_m$ has only fixed points,
Lemma~\ref{le:atr-repel} implies that for every $t\in\Gamma$, the
sequence $\alpha_{mn}(t)$ converges to a fixed point of $\alpha_m$
as $n\to\pm\infty$. Passing to the limit in (\ref{eq:lim-2}) as
$n\to\pm\infty$, in view of the continuity of $f$, we get
\begin{equation}
\left\{ \prod_{i=0}^{k-1} l_i \right\} \cdot \lim_{n\to\pm\infty}
 f_m[\alpha_{mn+k}(t)]
=
\left \{\prod_{i=0}^{k-1} l_i \right\} \cdot \lim_{n\to\pm\infty}
f_m[\alpha_{mn}(t)]  \label{eq:lim-3}
\end{equation}
where
\[
l_i:=f\left( \lim_{n\to\pm\infty}\alpha_{mn+i}(t) \right).
\]
Moreover, since the factors belonging to the first product from
the left (from the right) in (\ref{eq:lim-2}) are contained in the
second product from the right (from the left) in (\ref{eq:lim-2}),
we see that if $l_i=0$ for some $i=0,\dots,k-1$, then
\[
\lim_{n\to\pm\infty}  f_m[\alpha_{mn+k}(t)]
=
\lim_{n\to\pm\infty}
f_m[\alpha_{mn}(t)]
=0.
\]
Otherwise, we divide equality (\ref{eq:lim-3}) on the
product $l_0l_1\cdots l_{k-1}$, which again gives
(\ref{eq:lim-1}). \rule{2mm}{2mm}

For the  functional operator (\ref{eq:operator}), in the case $m \ge
1$ define the two continuous functions $\eta_i:\Gamma\to\bR, i\in \{
0,1 \}$, by the formulas
\begin{eqnarray}
\label{eq:eta0m}
\eta_0(t)&:=&|a_m(t)|-|b_m(t)|
\min\Big\{|\alpha_m'(t)|^{-\alpha_X},|\alpha_m'(t)|^{-\beta_X}\Big\},\\
\label{eq:eta1m}
\eta_1(t)&:=&|a_m(t)|-|b_m(t)|
\max\Big\{|\alpha_m'(t)|^{-\alpha_X},|\alpha_m'(t)|^{-\beta_X}\Big\}.
\end{eqnarray}
Clearly, the shift $\alpha_m$ has only fixed points. Then for
$m\ge 1$, we can define the following functions in correspondence
with Lemma~\ref{le:atr-repel}:
\begin{equation}\label{eq:etapmm}
\eta_i^\pm(t):=\lim_{n\to\pm\infty}\eta_i[\alpha_{mn}(t)],
\quad t\in\Gamma, \quad i\in\{0,1\}.
\end{equation}

In view of Lemma~\ref{le:atr-repel} and the continuity of the
functions $a_m, b_m, \alpha_m'$, the limits in (\ref{eq:etapmm})
exist for every $t\in\Gamma$. Hence, the functions $\eta_i^\pm$
are well-defined.
Introduce the sets $\Gamma_j, j\in\{1,2,\dots,5\}$,
as before, by formulas (\ref{eq:Gamma-def1}),
(\ref{eq:Gamma-def2}).
\begin{corollary}$\!\!\!${\bf .}\label{co:inv}
For every $i\in \{ 0,1\}$ and every $k\in\bZ$,
\[
\eta_i^-(t)=\eta_i^-[\alpha_k(t)], \quad
\eta_i^+(t)=\eta_i^+[\alpha_k(t)], \quad t\in\Gamma,
\]
and hence, the sets $\Gamma_j,\: j\in \{1,2,\dots, 5\}$, are
$\alpha_k$-invariant.
\end{corollary}

This statement follows from Proposition~\ref{pr:lim}.

Consider the function
\begin{equation}\label{eq:sigmam-def}
\sigma_A(t):=\left\{
\begin{array}{cl}
a_m(t)-b_m(t), & t\in\Gamma_1,\\
a_m(t),      & t\in\Gamma_2,\\
-b_m(t),     & t\in\Gamma_3,\\
0,         & t\in\Gamma\setminus(\Gamma_1\cup\Gamma_2\cup\Gamma_3).
\end{array}
\right.
\end{equation}
Obviously, the functions (\ref{eq:eta0m})--(\ref{eq:sigmam-def})
generalize the corresponding functions, introduced in Section 4,
to the case of arbitrary multiplicity $m\ge 1$ of periodic points
(recall that $m=2$ if $\alpha$ changes the orientation on $\Gamma$).

Now we are able to formulate the main result.
\begin{theorem}$\!\!\!${\bf .}\label{th:periodic}
Let $\Gamma$ be an oriented Jordan smooth curve and let $\alpha$
be a diffeomorphism of $\Gamma$ onto itself which has an arbitrary
nonempty set  of periodic points. The operator $A$ is right (left)
invertible in a reflexive r.-i. space $X(\Gamma)$ of fundamental
type with nontrivial Boyd indices if and only if
\begin{equation}\label{eq:periodic-right}
\sigma_A(t)\ne 0\mbox{ for all }
t\in\Gamma\setminus\Gamma_4\quad\mbox{ and }\quad R(\Gamma_4)\mbox{ holds}
\end{equation}
(respectively,
\begin{equation}\label{eq:periodic-left}
\sigma_A(t)\ne 0\mbox{ for all }
t\in\Gamma\setminus\Gamma_5\quad\mbox{ and }\quad L(\Gamma_5)\mbox{ holds )}.
\end{equation}
\end{theorem}

For the proof of this theorem we need some auxiliary results, which
will be stated in the next two subsections.
\subsection{Adjoint operator}\label{sec:adjoint}
Since the Lebesgue measure on $\Gamma$ is separable 
(see, e.g., \cite[Ch. 1, Subsection 6.10]{ka}),
from Lemma~\ref{le:ass&adj} it follows that the general form of a linear functional
on the reflexive r.-i. space $X(\Gamma)$ is given by
\begin{equation}\label{eq:form}
l(u)=(u,v):=\int_\Gamma u(\tau)\overline{v(\tau)}|d\tau|,
\quad u\in X(\Gamma),\quad v\in X'(\Gamma).
\end{equation}
From (\ref{eq:form}) it follows that the adjoint operator $A^*$
acting on the space $X'(\Gamma)=(X(\Gamma))^*$ and defined by the
equality $(Au,v)=(u,A^*v)$, has the form
\begin{equation}\label{eq:adjoint}
A^*=\overline{a}I-\overline{b(\alpha_{-1})}|\alpha_{-1}'|W^{-1}
\end{equation}
(here and in what follows $f(\alpha_k)=f\circ\alpha_k, \; k\in
{\bf Z}$).

It is easy to see that if $\tau$ is a periodic point of the
multiplicity $m\ge 1$ for the shift $\alpha$, then $\tau$ is also
a periodic point for the shift $\beta=\alpha_{-1}$ of the same
multiplicity.

Since the operator $A^*$ has the same form as the operator $A$,
for $A^*$ one can define functions $\widetilde{\eta}_i$ by analogy
with (\ref{eq:eta0m}) and (\ref{eq:eta1m}), and functions
$\widetilde{\eta}_i^\pm$ by analogy with (\ref{eq:etapmm}),
replacing the coefficients $a,\,b$, the shift $\alpha$, and the
Boyd indices $\alpha_X,\,\beta_X$, respectively, by the
coefficients $\overline{a}, \:
\overline{b(\alpha_{-1})}|\alpha_{-1}'|$, the shift $\alpha_{-1}$,
and the Boyd indices $\alpha_{X'},\,\beta_{X'}$ of the associate
(= dual) space $X'(\Gamma)$.
\begin{lemma}$\!\!\!${\bf .}\label{le:eta-adjoint}
For every $t\in\Gamma$ and $i\in \{0,1\}$, we have
\[
\widetilde{\eta}_i^+(t)=\eta_i^-(t),
\quad
\widetilde{\eta}_i^-(t)=\eta_i^+(t).
\]
\end{lemma}
{\bf Proof.} Taking into account the relations
\begin{equation}
\label{eq:eta-adjoint-1}
\prod_{i=0}^{m-1}\alpha'_{-1}[\alpha_{-i}(t)]=
\alpha_{-m}'(t)=(\alpha_{m}'[\alpha_{-m}(t)])^{-1}
\end{equation}
and the equalities (\ref{eq:ind-ass}), we infer from the
definition of $\widetilde{\eta}_0$ that
\begin{eqnarray}
\widetilde{\eta}_0(t)\! & = &\!  \left|
\prod_{i=0}^{m-1}\overline{a[\alpha_{-i}(t)]} \right|
-
\left| \prod_{i=0}^{m-1} \Big(
\overline{b[\alpha_{-i-1}(t)]}\cdot|\alpha_{-1}'[\alpha_{-i}(t)]|
\Big) \right|\min\Big\{ |\alpha_{-m}'(t)|^{-\alpha_{X'}},
|\alpha_{-m}'(t)|^{-\beta_{X'}} \Big\}  \nonumber\\  & = &\!
\left| a_m[\alpha_{-m+1}(t)] \right|
-
\left| b_m[\alpha_{-m}(t)] \right|\min\left\{
\left|\alpha'_{-m}(t)]\right|^{1-\alpha_{X'}},
\left|\alpha'_{-m}(t)]\right|^{1-\beta_{X'}} \right\} \nonumber\\
& = &\! \left| a_m[\alpha_{-m+1}(t)] \right|
-
\left| b_m[\alpha_{-m}(t)] \right|\min\left\{
\left|\alpha'_m[\alpha_{-m}(t)]\right|^{-\alpha_{X}},
\left|\alpha'_m[\alpha_{-m}(t)]\right|^{-\beta_{X}} \right\} .
 \label{eq:eta-adjoint-2}
\end{eqnarray}
Hence, from (\ref{eq:eta-adjoint-2}) we get
\begin{eqnarray}
\widetilde{\eta}_0[\alpha_{mn}(t)] &=& \left|
a_m[\alpha_{mn-m+1}(t)] \right| - \left| b_m[\alpha_{mn-m}(t)]
\right|\times \nonumber\\ &\times & \min\left\{
\left|\alpha_m'[\alpha_{mn-m}(t)]\right|^{-\alpha_{X}},
\left|\alpha_m'[\alpha_{mn-m}(t)]\right|^{-\beta_{X}} \right\}.
\label{eq:eta-adjoint-3}
\end{eqnarray}
Passing to the limit as $n\to+\infty$, we get from
(\ref{eq:eta-adjoint-3}) and Proposition~\ref{pr:lim} that
\[
\widetilde{\eta}_0^-(t)=
\lim_{n\to+\infty}\widetilde{\eta}_0[\alpha_{mn}(t)]
=
\lim_{n\to+\infty}\eta_0[\alpha_{mn}(t)]=\eta_0^+(t).
\]
Analogously, one can prove that
\[
\widetilde{\eta}_0^+(t)=\eta_0^-(t),
\quad
\widetilde{\eta}_1^-(t)=\eta_1^+(t),
\quad
\widetilde{\eta}_1^+(t)=\eta_1^-(t).
\quad
\rule{2mm}{2mm}
\]

Put $\widetilde{\Phi}:=\overline{\{t\in\Gamma:\alpha_{-m}(t)\ne
t\}}$. Clearly, $\widetilde{\Phi}=\Phi$. Consider the set
$\widetilde{\Gamma}_1:=\Gamma\setminus\widetilde{\Phi}=
\Gamma\setminus\Phi=\Gamma_1$. With the help of functions
$\widetilde{\eta}_i^\pm$, we define the sets
$\widetilde{\Gamma}_j,\; j\in \{2,3,4,5\}$, by analogy with
(\ref{eq:Gamma-def1}), (\ref{eq:Gamma-def2}). From
Lemma~\ref{le:eta-adjoint} we see that
\begin{equation}\label{eq:periodic-7}
\widetilde{\Gamma}_2=\Gamma_2,
\quad
\widetilde{\Gamma}_3=\Gamma_3,
\quad
\widetilde{\Gamma}_4=\Gamma_5,
\quad
\widetilde{\Gamma}_5=\Gamma_4.
\end{equation}
For the adjoint operator (\ref{eq:adjoint}), define the function
$\sigma_{A^*}(t)$ and the conditions $\widetilde{R}(\Omega),
\widetilde{L}(\Omega)$ by analogy with (\ref{eq:sigmam-def}) and
the conditions of $R(\Omega), L(\Omega)$, respectively, using the
coefficients $\overline{a}$ and
$\overline{b(\alpha_{-1})}|\alpha_{-1}'|$ instead of $a$ and $b$,
the shift $\alpha_{-1}$ instead of $\alpha$, the sets
$\widetilde{\Gamma}_j$ instead of $\Gamma_j$.

\begin{proposition}$\!\!\!${\bf .}\label{le:red2}
{\rm (a)}
The conditions
$L(\Gamma_5)$ and $\widetilde{R}(\widetilde{\Gamma}_4)$ are equivalent.

\noindent
{\rm (b)}
The conditions
$L(\Gamma_4)$ and $\widetilde{R}(\widetilde{\Gamma}_5)$ are equivalent.
\end{proposition}

The proof is straightforward in view of (\ref{eq:periodic-7}).
\begin{lemma}$\!\!\!${\bf .}\label{le:red1}
{\rm (a)}
$\sigma_A(t)\ne 0$ for all
$t\in\Gamma\setminus\Gamma_5$ if and only if $\sigma_{A^*}(t)\ne 0$
for all $t\in\Gamma\setminus\widetilde{\Gamma}_4$.

\noindent
{\rm (b)}
$\sigma_A(t)\ne 0$ for all
$t\in\Gamma\setminus\Gamma_4$ if and only if $\sigma_{A^*}(t)\ne 0$
for all $t\in\Gamma\setminus\widetilde{\Gamma}_5$.
\end{lemma}
{\bf Proof.} (a)
Since $\alpha_m(t)=t$ for $t\in\Gamma_1$, we have from the first identity
in (\ref{eq:eta-adjoint-1}) that
\[
\prod_{i=0}^{m-1}\alpha'_{-1}[\alpha_{-i}(t)]=\alpha_{-m}'(t)=1,
\quad t\in\widetilde{\Gamma}_1=\Gamma_1.
\]
Hence, using again the identity $\alpha_m(t)=t$ we get
\begin{eqnarray*}
\sigma_{A^*}(t)\! &=&\!
\prod_{i=0}^{m-1}\overline{a[\alpha_{-i}(t)]}
-
\prod_{i=0}^{m-1}\Big(\overline{b[\alpha_{-i-1}(t)]}\cdot|\alpha_{-1}'[\alpha_{-i}(t)]|\Big)
\!= \prod_{i=0}^{m-1}\overline{a[\alpha_{m-i}(t)]} -
\prod_{i=1}^{m}\overline{b[\alpha_{m-i}(t)]}
=
\\
\!&=&\!
\overline{a_m(t)-b_m(t)}
=
\overline{\sigma_A(t)}.
\end{eqnarray*}
Thus, $\sigma_{A^*}(t)\ne 0$ for $t\in\widetilde{\Gamma}_1$ if and
only if $\sigma_A(t)\ne 0$ for $t\in\Gamma_1$.

By Corollary~\ref{co:inv}, the set $\widetilde{\Gamma}_2=\Gamma_2$ is
$\alpha_{m-1}$-invariant. Consequently,
\[
\sigma_{A^*}(t)=\prod_{i=0}^{m-1}\overline{a[\alpha_{-i}(t)]}\ne 0,
\quad t\in\widetilde{\Gamma}_2,
\]
if and only if $\sigma_A(t)=a_m(t)\ne 0$ for $t\in\Gamma_2$.

Taking into account that $\alpha$ (and $\alpha_{-1}$) is
the diffeomorphism, we see that
\[
\sigma_{A^*}(t)=\prod_{i=0}^{m-1}
\Big(\overline{b[\alpha_{-i-1}(t)]}
\cdot|\alpha_{-1}'[\alpha_{-i}(t)]|\Big)\ne 0,
\quad t\in\widetilde{\Gamma}_3,
\]
if and only if
\begin{equation}\label{eq:periodic-8}
\prod_{i=1}^mb[\alpha_{-i}(t)]\ne 0,\quad t\in\widetilde{\Gamma}_3.
\end{equation}
Due to Corollary~\ref{co:inv}, the set $\widetilde{\Gamma}_3=\Gamma_3$ is
$\alpha_m$-invariant. Thus, (\ref{eq:periodic-8}) is equivalent to
$\sigma_A(t)=b_m(t)\ne 0$ for all $t\in\Gamma_3$.

Moreover, the condition $\sigma_{A^*}(t)\ne 0$ ($\sigma_A(t)\ne 0$) for all
$t\in\Gamma\setminus\widetilde{\Gamma}_4$ (respectively, $t\in\Gamma\setminus\Gamma_5$)
implies
$\Gamma=\widetilde{\Gamma}_1\cup\widetilde{\Gamma}_2
\cup\widetilde{\Gamma}_3\cup\widetilde{\Gamma}_4$
(respectively, $\Gamma={\Gamma}_1\cup{\Gamma}_2\cup{\Gamma}_3\cup{\Gamma}_5$).
In view of (\ref{eq:periodic-7}),
\[
\widetilde{\Gamma}_1\cup\widetilde{\Gamma}_2\cup\widetilde{\Gamma}_3
=
\Gamma\setminus\widetilde{\Gamma}_4
=
\Gamma\setminus{\Gamma}_5
=
{\Gamma}_1\cup{\Gamma}_2\cup{\Gamma}_3.
\]
Thus, $\sigma_A(t)\ne 0$ for all $t\in\Gamma\setminus\Gamma_5$
if and only if $\sigma_{A^*}(t)\ne 0$ for all $t\in\Gamma\setminus\widetilde{\Gamma}_4$.
The statement (b) is proved by analogy.
\rule{2mm}{2mm}
\subsection{Reduction to the case of fixed points}
In this subsection we reduce the investigation of the right
invertibility of the operator $A$ with the shift $\alpha$
preserving orientation on $\Gamma$ and having periodic points of
multiplicity $m>1$, or changing orientation on $\Gamma$ (then $m=2$), to the
investigation of the right invertibility of some operator of the
same form with the shift $\alpha_m$ which has only fixed points.
\begin{theorem}$\!\!\!${\bf .}\label{th:reduction}
The operator $A=aI-bW$ is right invertible in $X(\Gamma)$ if and
only if the operator $A_m:=a_mI-b_m(\alpha_{m-1})W^m$ is right
invertible in $X(\Gamma)$, and
\begin{equation}\label{eq:reduction-1}
|a_i(t)|+|b[\alpha_{i-1}(t)]|>0 \quad\mbox{for all}\quad
t\in\Gamma_4\quad\mbox{and all}\quad i=1,2,\dots,m-1.
\end{equation}
\end{theorem}
{\bf Proof.} If $m=1$, then the statement is trivial. So, assume
that $m>1$.

If the shift $\alpha$ preserves the orientation on $\Gamma$, then fix
some $\tau\in\Lambda$ and choose an arbitrary open arc $l$ from
the set
\[
\Big\{\Big(\tau,\alpha(\tau)\Big),\Big(\alpha(\tau),\alpha_2(\tau)\Big),
\dots,\Big(\alpha_{m-1}(\tau),\alpha_{m}(\tau)\Big)\Big\}.
\]
If the shift $\alpha$ changes the orientation on $\Gamma$, then
$\alpha$ has exactly two fixed points on $\Gamma$, say $z_1$ and
$z_2$. In this case take $m=2$, and let $l$ be an arbitrary open
arc from the set $\Big\{(z_1,z_2),(z_2,z_1)\Big\}$. Then
$\Gamma\setminus \{z_1,z_2 \}=l\cup\alpha(l)$.

Let  $[X(l)]^m$ be the space of all vectors with $m$ components
from $X(l)$ and $\sigma_l$ be the isomorphism of $X(\Gamma)$ onto
$[X(l)]^m$ defined by the rule
\[
(\sigma_l\varphi)(t):=\Big\{\varphi[\alpha_k(t)]\Big\}_{k=0}^{m-1},\quad t\in l.
\]
For $1\le i\le m-1$, define the operators $F_i:
[X(\Gamma)]^{m-i+1}\to[X(\Gamma)]^{m-i+1}$ by
\[
F_i:=\left[
\begin{array}{cccccc}
a_iI   & -b(\alpha_{i-1})I& \ddots & O                 & O \\
O      & a(\alpha_i)I     & \ddots & O                 & O \\
\ddots & \ddots           & \ddots & \ddots            & \ddots\\
O      & O                & \ddots & a(\alpha_{m-2})I  & -b(\alpha_{m-2})I \\
-b_i(\alpha_{m-1})W^m & O & \ddots & O                 & a(\alpha_{m-1})I
\end{array}
\right].
\]
Since
\[
\sigma_l A\sigma_l^{-1}\varphi= F_1\varphi\quad\mbox{for
all}\quad\varphi\in[X(l)]^m,
\]
and since $l$ is an arbitrary connected component in
$\Gamma\setminus \{ \tau,\alpha(\tau),\dots,\alpha_{m-1}(\tau)\}$
or $\Gamma\setminus \{ z_1,z_2 \}$, the operator $A$ is right
invertible in $X(\Gamma)$ if and only if the operator $F_1$ is
right invertible in each space
$[X(\alpha_k(l))]^m,k\in\{0,1,\dots,m-1\}$. In view of the
equality
\[
[X(\Gamma)]^m=
[X(l)]^m\oplus[X(\alpha(l))]^m
\oplus\dots\oplus[X(\alpha_{m-1}(l))]^m,
\]
the operator $A$ is right invertible in the space $X(\Gamma)$ if
and only if the operator $F_1$ is right invertible in the space
$[X(\Gamma)]^m$. If the operator $F_i,\; i\in\{1,2,\dots,m-1\}$,
is right invertible in the space $[X(\Gamma)]^{m-i+1}$, then we
get the $i$-th relation in
\begin{equation}\label{eq:reduction-prime}
\min_{t\in\Gamma}\Big(|a_i(t)|+|b[\alpha_{i-1}(t)]|\Big)>0, \quad
i=1,2,\dots,m-1.
\end{equation}
Indeed, otherwise all the elements of the $i$-th row have zero at
a point $t\in\Gamma$, which implies that ${\rm Im}\, F_i\ne
[X(\Gamma)]^{m-i+1}$.

If (\ref{eq:reduction-prime}) holds for some $i \in \{ 1,2,
\dots,m-1\}$, then there is a function $g^{(i)}\in C(\Gamma)$ such
that the function $f^{(i)}:=a_ig^{(i)}+b(\alpha_{i-1})$ is bounded
away from zero. Consequently, the operator
\[
C_i:=\left[
\begin{array}{cc}
a_iI & -b(\alpha_{i-1})I \\[3mm] (1/f^{(i)})I & (g^{(i)}/f^{(i)})I
\end{array}
\right], \quad i=1,2,\dots,m-1,
\]
is invertible in the space $[X(\Gamma)]^2$, and the following
equality holds:
\begin{equation}\label{eq:reduction-3}
F_i=\left \{
  \begin{array}{lll}
  \left[
\begin{array}{cc}
I & O_{1\times(m-i)}\\
D_i & F_{i+1}
\end{array}
\right]
\left[
\begin{array}{cc}
C_i & O_{2\times(m-i-1)}\\
O_{(m-i-1)\times 2} & I_{m-i-1}
\end{array}
\right], & & i=1,2,\dots,m-2,\\[5mm] 
\left[
\begin{array}{cc}
I & O\\ D_{i} & F_{i+1}
\end{array}
\right] C_{i},& & i=m-1,
\end{array} \right.
\end{equation}
where the operator $D_i:X(\Gamma)\to[X(\Gamma)]^{m-i}$ is defined by the formula
\[
D_i:=
\left\{
\begin{array}{ll}
{\rm column}\Biggl[-\frac{\displaystyle a(\alpha_i)}{\displaystyle
f^{(i)}}I,\overbrace{O,\dots,O}^{m-2-i},
-b_i(\alpha_{m-1})W^m\frac{\displaystyle g^{(i)}}{\displaystyle
f^{(i)}}I\Biggr], & i=1,2,\dots,m-2,
\\ 
-\frac{\displaystyle a(\alpha_{m-1})}
      {\displaystyle f^{(m-1)}}I-b_{m-1}(\alpha_{m-1})W^m
\frac{\displaystyle g^{(m-1)}}{\displaystyle f^{(m-1)}}I, 
&
i=m-1,
\end{array}
\right.
\]
$I_k$ is the identity operator in the space $[X(\Gamma)]^k$, $O_{k\times p}$ is the
zero operator from the space $[X(\Gamma)]^p$ into the space $[X(\Gamma)]^k$, and $F_m:=A_m$.

Since the second term on the right of (\ref{eq:reduction-3}) is
invertible in the space $[X(\Gamma)]^{m-i+1}$, the right
invertibility of the operator $F_i$ in the space
$[X(\Gamma)]^{m-i+1}$ is equivalent to the right invertibility of
the first term on the right of (\ref{eq:reduction-3}) in the space
$[X(\Gamma)]^{m-i+1}$, which in its turn is equivalent to the
right invertibility of the operator $F_{i+1}$ in the space
$[X(\Gamma)]^{m-i}$.

Indeed, if $F_{i+1}$ is right invertible in $[X(\Gamma)]^{m-i}$
and $F_{i+1}^{(-1)}$ is one of its right inverses, then
\[
\left[
\begin{array}{cc}
I   & O_{1\times(m-i)}\\
D_i & F_{i+1}
\end{array}
\right]
\left[
\begin{array}{cc}
I   & O_{1\times(m-i)}\\
-F_{i+1}^{(-1)}D_i & F_{i+1}^{(-1)}
\end{array}
\right]
=
\left[
\begin{array}{cc}
I   & O_{1\times(m-i)}\\
O_{(m-i)\times 1} & I_{m-i}
\end{array}
\right].
\]
Hence, the first term on the right of (\ref{eq:reduction-3}) is
right invertible.
On the other hand, if
\begin{eqnarray*}
\left[\!
\begin{array}{cc}
I   & O_{1\times(m-i)}\\
D_i & F_{i+1}
\end{array}\!
\right] \left[
\begin{array}{cc}
B_1 & B_2\\
B_3 & B_4
\end{array}
\right]\!\! &=&\! \!\left[\!
\begin{array}{cc}
B_1 & B_2\\
D_iB_1+F_{i+1}B_3 & D_iB_2+F_{i+1}B_4
\end{array} \!
\right]
\left[\!
\begin{array}{cc}
I   & O_{1\times(m-i)}\\
O_{(m-i)\times 1} & I_{m-i}
\end{array}\!
\right],
\end{eqnarray*}
then $B_2=O_{1\times(m-i)}$. Hence, $F_{i+1}B_4=I_{m-i}$, that is,
the operator $F_{i+1}$ is right invertible in the space
$[X(\Gamma)]^{m-i}$.

Thus, the operator $A$ is right invertible in the space
$X(\Gamma)$ if and only if the operator $A_m=F_m$ is right
invertible in the space $X(\Gamma)$ and, for every
$i\in\{1,2,\dots,m-1\}$, inequality (\ref{eq:reduction-prime})
holds. Clearly, (\ref{eq:reduction-prime}) implies
(\ref{eq:reduction-1}), so the necessity is proved.

It remains to prove that (\ref{eq:reduction-1}) and the right
invertibility of $A_m$ imply (\ref{eq:reduction-prime}). From
(\ref{eq:eta0m})--(\ref{eq:etapmm}),
(\ref{eq:eta0})--(\ref{eq:eta1}), (\ref{eq:etapm-def}) and
Proposition~\ref{pr:lim} it follows that the functions
$\eta_i^\pm,i\in \{0,1\}$, defined for operators
$A_m=a_mI-b_m(\alpha_{m-1})W_{\alpha_m}$ and $A$, coincide. Hence,
for these operators, the corresponding sets $\Gamma_j,\; j\in
\{1,2,\dots,5\}$, coincide as well.

If the operator $A_m$ is right invertible, then, by Theorem~\ref{th:fixed2},
\begin{equation}\label{eq:reduction-4}
\sigma_{A_m}(t)\ne 0
\quad
\mbox{for all}\quad t\in\Gamma\setminus\Gamma_4,
\end{equation}
where
\[
\sigma_{A_m}(t):=\left\{
\begin{array}{cl}
a_m(t)-b_m[\alpha_{m-1}(t)], & t\in\Gamma_1,\\
a_m(t),      & t\in\Gamma_2,\\
-b_m[\alpha_{m-1}(t)],     & t\in\Gamma_3,\\
0,         & t\in\Gamma\setminus(\Gamma_1\cup\Gamma_2\cup\Gamma_3).
\end{array}
\right.
\]
It is easy to see that $b_m[\alpha_{m-1}(t)]=b_m(t)$ for
$t\in\Gamma_1$. In view of Corollary~\ref{co:inv}, the set
$\Gamma_3$ is $\alpha_{m-1}$-invariant. Hence,
$b_m[\alpha_{m-1}(t)]\ne 0$ for all $t\in\Gamma_3$ if and only if
$b_m(t)\ne 0$ for all $t\in\Gamma_3$. Thus, (\ref{eq:reduction-4})
is equivalent to
\begin{equation}\label{eq:reduction-5}
\sigma_{A}(t)\ne 0
\quad
\mbox{for all}\quad t\in\Gamma\setminus\Gamma_4.
\end{equation}
From (\ref{eq:reduction-5}) it follows that
\begin{equation}\label{eq:reduction-6}
|a_i(t)|+|b[\alpha_{i-1}(t)]|>0\quad\mbox{for all} \quad
t\in\Gamma\setminus\Gamma_4\quad\mbox{and all}\quad
i=1,2,\dots,m-1.
\end{equation}
Indeed, assume the contrary: for some $i\in\{1,2,\dots,m-1\}$
there is a point
 $t_0\in\Gamma\setminus\Gamma_4$ such that $|a_i(t_0)|=|b[\alpha_{i-1}(t_0)]|=0$.
Hence, $a_m(t_0)=0$ and $b_m(t_0)=0$. Consequently,
$\sigma_A(t_0)=0$ and we get the contradiction. Clearly,
(\ref{eq:reduction-1}) and (\ref{eq:reduction-6}) imply
(\ref{eq:reduction-prime}), which completes the proof of the
sufficiency. \rule{2mm}{2mm}
\subsection{Proof of Theorem~\ref{th:periodic}}
The case of $m=1$ was considered in Theorem~\ref{th:fixed2}.
Assume that $m>1$ and consider the case of the right invertiblity
of $A$.

In view of Theorem~\ref{th:reduction}, the right invertibility of
$A$ is equivalent to (\ref{eq:reduction-1}) and the right
invertibility of $A_m$. By Theorem~\ref{th:fixed2}, the right invertibility
of $A_m$ is equivalent to (\ref{eq:reduction-4}) and the condition
\begin{equation}\label{eq:periodic-1}
\begin{array}{cc}
\mbox{for every }t\in\Gamma_4\mbox{ there is }k_0\in\bZ\mbox{ such
that }\\ a_m[\alpha_{mk}(t)]\ne 0\mbox{ for }k<k_0, \quad
b_m[\alpha_{mk+m-1}(t)]\ne 0\mbox{ for }k\ge k_0.
\end{array}
\end{equation}
Let us show that the conjunction of (\ref{eq:reduction-1}) and
(\ref{eq:periodic-1}) is equivalent to the condition
\begin{equation}\label{eq:periodic-2}
\begin{array}{cc}
\mbox{for every }t\in\Gamma_4\mbox{ there are }
s\in\{0,1,\dots,m-1\}\mbox{ and } k_0\in\bZ\mbox{ such that }\\
a[\alpha_{k}(t)]\ne 0\mbox{ for }k<mk_0+s, \quad
b[\alpha_{k}(t)]\ne 0\mbox{ for }k\ge mk_0+s.
\end{array}
\end{equation}
Indeed, condition (\ref{eq:periodic-1}) is equivalent to the condition
\begin{equation}\label{eq:periodic-3}
\begin{array}{cc}
\mbox{for every }t\in\Gamma_4\mbox{ there is }k_0\in\bZ\mbox{ such that }\\
a[\alpha_{k}(t)]\ne 0\mbox{ for }k<mk_0,
\quad
b[\alpha_{k}(t)]\ne 0\mbox{ for }k\ge mk_0+m-1.
\end{array}
\end{equation}
Since the set $\Gamma_4$ is invariant with respect to the shift $\alpha_{mk_0}$,
(\ref{eq:reduction-1}) is equivalent to the condition
\begin{equation}\label{eq:periodic-4}
\begin{array}{c}
\mbox{for every }t\in\Gamma_4\mbox{ and every }
i\in\{1,2,\dots,m-1\},\\[1mm] \prod_{j=0}^{i-1}
 |a[\alpha_{mk_0+j}(t)]|
+|b[\alpha_{mk_0+i-1}(t)]|>0.
\end{array}
\end{equation}
Fix $t\in\Gamma_4$. If $a[\alpha_{mk_0+i}(t)]\ne 0$ for all $i \in
\{ 0,1,\dots, m-2 \}$, then from (\ref{eq:periodic-3}) we get
(\ref{eq:periodic-2}) with $s=m-1$. Otherwise there exists an
$s\in\{0,1,\dots,m-2\}$ such that
\begin{equation}\label{eq:periodic-5}
a[\alpha_{mk_0+i}(t)]\ne 0,\quad i\in\{0,1,\dots,s-1\},\quad
a[\alpha_{mk_0+s}(t)]=0.
\end{equation}
Then $\prod_{j=0}^{i-1}a[\alpha_{mk_0+j}(t)]=0$ for
$i\in\{s+1,\dots,m-1\}$. In that case (\ref{eq:periodic-4})
implies
\begin{equation}\label{eq:periodic-6}
b[\alpha_{mk_0+i-1}(t)]\ne 0,\quad i\in\{s+1,\dots,m-1\}.
\end{equation}
From (\ref{eq:periodic-5}), (\ref{eq:periodic-6}) and
(\ref{eq:periodic-3}) we get (\ref{eq:periodic-2}) with $s\in
\{0,1,\ldots,m-2\}$.

On the other hand, if (\ref{eq:periodic-2}) holds, then the
conditions (\ref{eq:periodic-3}) and (\ref{eq:periodic-4}) are
fulfilled. But as it was said above, (\ref{eq:periodic-4}) is
equivalent to (\ref{eq:reduction-1}), and (\ref{eq:periodic-3}) is
equivalent to (\ref{eq:periodic-1}). So, we get that we need. It
remains to note that condition (\ref{eq:periodic-2}) is equivalent
to $R(\Gamma_4)$.

Thus, we have proved that the right invertibility of $A$ is
equivalent to (\ref{eq:reduction-4}) and the property
$R(\Gamma_4)$. But, as was shown in the proof of
Theorem~\ref{th:reduction}, the properties (\ref{eq:reduction-4})
and (\ref{eq:reduction-1}) are equivalent. The case of the right
invertibility is considered.

The case of the left invertibility is reduced to the previous one
by passing to adjoint operators and applying
Proposition~\ref{le:red2}(a) and Lemma~\ref{le:red1}(a).
\rule{2mm}{2mm}
\begin{corollary}$\!\!\!${\bf .}\label{co:invertibility}
The operator $A$ is invertible in 
$X(\Gamma)$ if and only if
$\sigma_A(t)\ne 0$ for all $t\in\Gamma$.
\end{corollary}
{\bf Proof.} The condition $\sigma_A(t)\ne 0$ for all
$t\in\Gamma\setminus\Gamma_4$ (respectively,
$t\in\Gamma\setminus\Gamma_5$) implies that $\Gamma_5=\emptyset$
(respectively, $\Gamma_4=\emptyset$). Thus, in the case of the
invertibility of $A$, from Theorem~\ref{th:periodic} we get
$\Gamma_4=\Gamma_5=\emptyset$, and conditions
(\ref{eq:periodic-left}), (\ref{eq:periodic-right}) degenerate to
$\sigma_A(t)\ne 0$ for all
$t\in\Gamma=\Gamma_1\cup\Gamma_2\cup\Gamma_3$.

On the other hand, the condition $\sigma_A(t)\ne 0$ for all $t\in\Gamma$
implies that $\Gamma=\Gamma_1\cup\Gamma_2\cup\Gamma_3$, that is,
$\Gamma_4=\Gamma_5=\emptyset$. Hence, by Theorem~\ref{th:periodic},
the operator $A$ is simultaneously right and left invertible.
\rule{2mm}{2mm}
\subsection{Spectrum of the weighted shift operator}\label{sec:spectrum}
The spectrum of the weighted shift operator is calculated in
\cite[Lemma~6.10]{AKarl00} in the case of only two fixed points.
Now we generalize this result to the case of arbitrary nonempty set of
periodic points.

For a continuous function $d$ and a number $m\in\bN$, consider the
two functions
\begin{eqnarray*}
\delta(t):=
|d_m(t)|\min\Big\{|\alpha_m'(t)|^{-\alpha_X},|\alpha_m'(t)|^{-\beta_X}\Big\},
\\
\Delta(t):=
|d_m(t)|\max\Big\{|\alpha_m'(t)|^{-\alpha_X},|\alpha_m'(t)|^{-\beta_X}\Big\}.
\end{eqnarray*}
\begin{theorem}$\!\!\!${\bf .}\label{th:spectrum}
The spectrum of the weighted shift operator $dW$ with a weight (coefficient)
$d\in C(\Gamma)$ in a reflexive r.-i. space $X(\Gamma)$ of fundamental type with nontrivial
Boyd indices $\alpha_X,\beta_X$ has the form
\[
\sigma(dW)= \left( \bigcup_{t\in\Gamma\setminus\Phi} \Big\{
z\in\bC:z^m=d_m(t) \Big\} \right)\cup \left(
\bigcup_{\gamma\subset\Phi\setminus\Lambda}\Omega(d,\gamma)
\right) \cup \left( \bigcup_{\tau\in Y'}\Psi(d,\tau) \right),
\]
where $\gamma$ are connected components of $\Phi\setminus
\Lambda$,
\[
\Omega(d,\gamma) =
\left\{
\begin{array}{lll}
\Big\{
z\in\bC:
\min\limits_{\tau\in\Lambda\cap\overline{\gamma}}\delta(\tau)
\le |z|^m\le
\max\limits_{\tau\in\Lambda\cap\overline{\gamma}}\Delta(\tau)
\Big\}
& \mbox{if} & d_m\in GC(\overline{\gamma}),
\\
\Big\{
z\in\bC:
 |z|^m\le
\max\limits_{\tau\in\Lambda\cap\overline{\gamma}}\Delta(\tau)
\Big\}
& \mbox{if} & d_m\not\in GC(\overline{\gamma}),
\end{array}
\right.
\]
\[
\Psi(d,\tau)=\Big\{
z\in\bC:
\delta(\tau)\le |z|^m\le\Delta(\tau)
\Big\}.
\]
\end{theorem}
{\bf Proof.} Corollary~\ref{co:invertibility} and the equivalence
of (\ref{eq:reduction-5}) and (\ref{eq:reduction-4}) imply that
the operator $B:=zI-dW, z\in\bC$, is not invertible in $X(\Gamma)$ if and only if
the operator $B_m=z^mI-d_m W_{\alpha_m}$ is not invertible in $X(\Gamma)$, where
$m$ is the multiplicity of periodic points of the shift $\alpha$
if $\alpha$ preserves the orientation on $\Gamma$, and $m=2$
otherwise. Clearly, the shift $\alpha_m$ has only fixed points,
and we can apply Theorem~\ref{th:fixed} to the operator $B_m$. It
is easy to see that
\[
\eta_0(t)=|z|^m-\delta(t),\quad \eta_1(t)=|z|^m-\Delta(t).
\]
Hence, by Theorem~\ref{th:fixed}, the spectrum of the operator $dW$ has the form
\[
\begin{array}{l}
\sigma(dW) = \left( \bigcup_{\gamma\subset\Phi\setminus\Lambda}
\Big\{ z\in\bC:z^mI-d_mW_{\alpha_m}\quad\mbox{is not invertible
in}\quad X(\gamma) \Big\} \right)
\\[4mm]\cup\left(
\bigcup_{t\in\Gamma\setminus\Phi}\Big\{z\in\bC:z^m=d_m(t)\Big\}\right)
 \cup\left( \bigcup_{\tau\in Y'} \Big\{ z\in\bC:
(|z|^m-\delta(\tau))(|z|^m-\Delta(\tau))\le 0 \Big\} \right).
\end{array}
\]
By \cite[Lemma~6.10]{AKarl00} and Lemma~\ref{le:subset}, for
$\gamma \in \Phi\setminus \Lambda$,
\[
 \Big\{ z\in\bC:z^mI-d_mW_{\alpha_m}\quad\mbox{is
not invertible in}\quad X(\gamma) \Big\}=\Omega(d,\gamma).
\]
Since $0\le\delta(t)\le\Delta(t)$, the inequality
$(|z|^m-\delta(t))(|z|^m-\Delta(t))\le 0$ is equivalent to the
inequality $\delta(t)\le |z|^m\le\Delta(t)$. Hence, for $\tau \in
Y'$,
\[
 \Big\{ z\in\bC:
(|z|^m-\delta(\tau))(|z|^m-\Delta(\tau))\le 0
\Big\}=\Psi(d,\tau).\quad \rule{2mm}{2mm}
\]

In particular, if $\gamma=(\tau_-,\tau_+)$ is an arc, the
endpoints of which are fixed points of the shift $\alpha$ and
$\Lambda\cap\gamma=\emptyset$, then from the description of the
set $\Omega(1,\gamma)$ we see that the spectral radius of the
shift operator $W$ in $X(\gamma)$ is given by (see also
\cite[Lemma~6.10]{AKarl00})
\begin{eqnarray*}
r(W;X(\gamma))
&=&\max_{\tau\in\{\tau_-,\tau_+\}}\max\{|\alpha'(\tau)|^{-\alpha_X},|\alpha'(\tau)|^{-\beta_X}\}
\\ &=&
\max\Big\{r(W;L^{1/\alpha_X}(\gamma)),r(W;L^{1/\beta_X}(\gamma))\Big\}.
\end{eqnarray*}
This simple example shows, that estimate (\ref{eq:spectral-est}) in Theorem~\ref{th:spectral-est}
is sharp.

Note that Theorem~\ref{th:spectrum} were obtained for the case of reflexive Orlicz spaces by
V.~Aslanov and the second author \cite{AK89} (see also \cite{AK2000}, \cite[Corollary~2.1]{Karl95}).

\vspace{2mm} \noindent {\bf Remark.} We do not use the restriction
on r.-i. spaces to be of fundamental type in the proofs of all
results of this paper, except for Theorem~\ref{th:two-fixed}. But
the proof of Theorem~\ref{th:two-fixed} (see \cite{AKarl00})
essentially depends on this restriction.

\vspace{2mm}
\noindent
{\bf Problem.}
Is it possible to drop the restrictions $\alpha_X=p_X$ and $\beta_X=q_X$?
\vskip 1truecm

\baselineskip=12pt


\vspace{1cm}
\begin{minipage}[t]{7cm}
A.~Yu.~Karlovich\\
Departamento de Matem\'atica\\
Instituto Superior T\'ecnico\\
Av. Rovisco Pais\\
1049-001, Lisboa,\\
Portugal\\
\\
E-mail: akarlov@math.ist.utl.pt\\
http://www.math.ist.utl.pt/$\widetilde{ }$akarlov
\end{minipage}
\hspace{8mm}
\begin{minipage}[t]{8cm}
Yu.~I.~ Karlovich\\
CINVESTAV del I.P.N.\\
Departamento de Matem{\'a}ticas\\
Apartado Postal 14-740\\
07000, Mexico, D.F., \\
Mexico\\
\\
E-mail: karlovic@math.cinvestav.mx
\end{minipage}

\vspace{1cm}
MSC 1991: Primary 39B32, 47B38

\hspace{2cm}
Secondary 46E30, 47A10

\begin{thebibliography}{10}
\bibitem{AAK}
{\sc Abramovich,~Yu., Arenson,~E.}, and {\sc Kitover,~A.:} {\it
Banach $C(K)$-modules and operators preserving disjointness}.
Pitman Research Notes in Mathematics Series {\bf 277}, Longman
Scientific \& Technical, New York 1992.

\bibitem{Antonevich}
{\sc Antonevich,~A.~B.:} {\it Linear Functional Equations.
Operator Approach}. University Press, Minsk 1988 (Russian).
English transl.: Birkh\"auser Verlag, Basel, Boston, Berlin 1995.

\bibitem{AntLeb}
{\sc Antonevich, A.,} and {\sc Lebedev, A.:} {\it
Functional-Differential Equations: {\rm {I}}. ${C}\sp *$-theory}.
Longman Scientific \& Technical, Harlow 1994.

\bibitem{AntLebBel}
{\sc Antonevich,~A., Belousov,~M.,}  and {\sc Lebedev,~A.:} 
{\it Functional Differential Equations: {\rm {I}{I}}.
${C}\sp*$-applications. {P}arts {\rm 1--2}}. 
Longman, Harlow 1998.

\bibitem{AK2000}
{\sc Aslanov,~V.:} Functional and Singular Integral Operators with
a Shift in Orlicz Spaces. Ph. D. Thesis, Baku, 1992 (Russian).

\bibitem{AK89}
{\sc Aslanov,~V.}, and {\sc Karlovich,~Yu.~I.:} One-sided
invertibility of functional operators in reflexive Orlicz spaces.
{\it Dokl. Akad. Nauk AzSSR} {\bf 45} (1989), no. 11-- 12, 3--7
(Russian).

\bibitem{BeSh}
{\sc Bennett,~C.}, and {\sc Sharpley,~R.:} {\it Interpolation of
Operators}. Academic Press, Boston 1988.

\bibitem{b4}
{\sc Boyd,~D.~W.:} Indices of function spaces and their
relationship to interpolation. {\it Canad. J. Math.} {\bf 21}
(1969), 1245--1254.

\bibitem{Lat}
{\sc Chicone,~C.,} and {\sc Latushkin,~Yu.:} {\it Evolution
Semigroups in Dynamical Systems and Differential Equations}.
American Mathematical Society, Providence, RI, 1999.

\bibitem{feh83}
{\sc Feh\'er,~F.:} Indices of Banach function spaces and spaces of
fundamental type. {\it J. Approx. Theory} {\bf 37} (1983), 12--28.

\bibitem{gk}
{\sc Gohberg,~I.}, and {\sc Krupnik,~N.:} {\it One-Dimensional
Linear Singular Integral Equations}. Vol. 1. Birkh\"auser
Verlag, Basel, Boston, Berlin 1992 (Russian original: Shtiintsa,
Kishinev 1973).

\bibitem{ka}
{\sc Kantorovich,~L.~V.}, and {\sc Akilov,~G.~P.}: {\it Functional
analysis}. Nauka, Moscow, 3rd ed., 1984 (in Russian). English
transl.: Pergamon Press, Oxford, 2nd ed., 1982.

\bibitem{AKarl00}
{\sc Karlovich,~A.~Yu.:} Criteria for one-sided invertibility of a
functional operator in rearrangement-invariant spaces of
fundamental type, to appear in {\it Mathematische Nachrichten}.

\bibitem{Karl95}
{\sc Karlovich,~Yu.~I.:} The continuous invertibility of
functional operators in Banach spaces. {\it Dissertationes
Math. (Rozprawy Mat.)} {\bf 340} (1995), 115--136.

\bibitem{KK76}
{\sc Karlovich,~Yu.~I.}, and {\sc Kravchenko,~V.~G.:} A Noether
theory for a singular integral operator with a shift having
periodic points. {\it Soviet Math. Dokl.} {\bf 17} (1976),
1547--1551.

\bibitem{KrLit}
{\sc Kravchenko,~V.~G.}, and {\sc Litvinchuk,~G.~S.:} {\it
Introduction to the Theory of Singular Integral Operators with
Shift}. Series: Mathematics and its applications {\bf 289}, Kluwer
Academic Publishers, Dordrecht, Boston, London 1994.

\bibitem{kps}
{\sc Krein,~S.~G., Petunin,~Ju.~I.}, and {\sc Semenov,~E.~M.:}
{\it Interpolation of Linear Operators}. Nauka, Moscow 1978
(Russian). English transl.: AMS Translations of Mathematical
Monographs {\bf 54}, Providence, R.I., 1982.

\bibitem{Kurbatov}
{\sc Kurbatov,~V.~G.:} {\it Functional Differential Operators and
Equations}. Series: Mathematics and its applications {\bf 473},
Kluwer Academic Publishers, Dordrecht, Boston, London 1999.

\bibitem{LT}
{\sc Lindenstrauss,~J.} and {\sc Tzafriri,~L.:} {\it Classical
Banach Spaces. Function Spaces}. Springer Verlag, New York, Berlin
1979.

\bibitem{Lit}
{\sc Litvinchuk,~G.~S.:} {\it Boundary Value Problems and Singular
Integral Equations with Shift}. Nauka, Moscow 1977 (Russian).

\bibitem{mal}
{\sc Maligranda,~L.:} Indices and interpolation, {\it
Dissertationes Math. (Rozprawy Mat.)} {\bf 234} (1985), 1--49.

\bibitem{mal1}
{\sc Maligranda,~L.:} {\it Orlicz Spaces and Interpolation}. Sem.
Math. 5, Dep. Mat., Univ. Estadual de Campinas, Campinas SP,
Brazil, 1989.

\bibitem{Mar85}
{\sc Mardiev,~R.:} A criterion for the semi-Noetherian property of
one class of singular integral operators with a non-Carleman
shift. {\it Dokl. Akad Nauk UzSSR} (1985), no. 2,  5--7 (Russian).

\bibitem{Mar88}
{\sc Mardiev,~R.:} A criterion for $n(d)$-normality of singular
integral operators with a shift having periodic points in Lebesgue
spaces. Samarkand, 1988, 41 pp. Manuscript no. 821-Uz88, deposited
at UzNIINTI (Russian).

\bibitem{zippin}
{\sc Zippin,~M.:} Interpolation of operators of weak type between
rearrangement invariant spaces. {\it J. Functional Analysis} {\bf
7} (1971), 267--284.
\end{thebibliography}
\end{document}